\documentclass[11pt]{article}
\usepackage[margin=1in]{geometry}
\usepackage[utf8]{inputenc}
\usepackage{amsthm}
\usepackage[affil-it]{authblk}
\usepackage{amsmath,amsfonts}
\usepackage{algorithmic}
\usepackage{algorithm}
\usepackage{array}
\usepackage[caption=false]{subfig}
\usepackage{textcomp}
\usepackage{url}
\usepackage{verbatim}
\usepackage{graphicx}
\usepackage{cite}
\usepackage{booktabs}

\usepackage{pdflscape}
\usepackage{placeins}
\usepackage[table,dvipsnames]{xcolor}
\usepackage{amsfonts}
\usepackage{amsmath} 
\usepackage[T1]{fontenc}
\usepackage{adjustbox}

\usepackage{hyperref}
\usepackage{float}
\usepackage{comment}
\usepackage{xurl}

\usepackage{enumitem}
\setlist[itemize]{noitemsep, topsep=0pt}
\setlist[enumerate]{noitemsep, topsep=0pt}

\newtheorem{theorem}{Theorem}[section]

\newtheorem{proposition}[theorem]{Proposition}

\newtheorem{remark}[theorem]{Remark}

\newtheorem{definition}{Definition}
\newtheorem{example}[theorem]{Example}

\DeclareMathOperator {\V}{\mathbb{V}} 
\DeclareMathOperator {\E}{\mathbb{E}} 
\DeclareMathOperator {\Skew}{\mathbb{S}} 

\begin{document}

\title{Quantile-Based Skewness for Fuzzy Numbers With Probabilistic Foundations: An Application in Portfolio Optimization}

\author[1]{Jan Schneider}
\author[1]{Kaja Bilińska\thanks{Corresponding author. E-mail: \href{mailto:Kaja.Bilinska@pwr.edu.pl}{Kaja.Bilinska@pwr.edu.pl}}}
\author[1]{Tomasz Szandała}
\author[2]{Paul Schneider}
\affil[1]{Wrocław University of Science and Technology, 50-370 Wrocław, Poland\protect\\
\texttt{\{Jan.Schneider, Kaja.Bilinska, Tomasz.Szandala\}@pwr.edu.pl}}
\affil[2]{Università della Svizzera italiana, Via Buffi 13, 6900 Lugano, Switzerland\protect\\
\texttt{paul.schneider@usi.ch}}
\date{}

\newcommand{\ieeenotice}{\copyright~2026 IEEE. Personal use of this material is permitted.
  Permission from IEEE must be obtained for all other uses, in any current or future media,
  including reprinting/republishing this material for advertising or promotional purposes,
  creating new collective works, for resale or redistribution to servers or lists, or reuse
  of any copyrighted component of this work in other works.
  This is the accepted version of the article published in \emph{IEEE Transactions on Fuzzy
  Systems}, vol.~34, no.~10, pp.~3362--3373, October 2026,
  DOI: \href{https://doi.org/10.1109/TFUZZ.2026.3721516}{10.1109/TFUZZ.2026.3721516}.}

\maketitle
\begingroup\renewcommand\thefootnote{}\footnotetext{\ieeenotice}\addtocounter{footnote}{-1}\endgroup

\begin{abstract}

This article quantifies asymmetry under imprecision by introducing a probabilistic interpretation for fuzzy numbers and a nonmoment-based skewness coefficient derived from it. The coefficient adapts classical quantile-based skewness measures (Bowley, Kelley, and Groeneveld--Meeden) to the $\alpha$-cut representation of fuzzy numbers. Existing studies typically employ heuristic modifications of the classical Fisher--Pearson formula, therein substituting membership functions for probability densities without a rigorous semantic bridge between the two. Our coefficient, in contrast, is grounded in am explicit probabilistic interpretation, rendering it both mathematically meaningful and aligned with the semantics of fuzzy set theory.

Our approach interprets a fuzzy number's left and right membership function components as cumulative and survival probability functions of associated random variables. This provides a probabilistic reading of its $\alpha$-cuts as generalized quantiles representing intervals of values that are ``at least $(1-\alpha)$-probable,'' thereby highlighting a duality between probability and possibility theory.
As a quantile-based measure, the proposed coefficient offers invariance under scale and location transformations. Moreover, its point-based variant requires only a fixed number of quantile evaluations per computation, constant time in the $\alpha$-discretization, compared with the linear cost of moment-based alternatives.
The coefficient comprises two complementary constituents: an ``inner'' measure quantifying the intrinsic skewness of the underlying probabilistic distributions, and an ``outer'' measure capturing the fuzzy number's overall profile asymmetry. We illustrate the measures within a fuzzy mean-variance-skewness portfolio optimization framework on historical equity data, comparing their performance with two of the most highly cited original moment-based fuzzy skewness coefficients from the literature.

\end{abstract}

\medskip
\noindent\textbf{Index Terms:} Financial engineering, fuzzy numbers, Markowitz portfolio optimization, mathematical finance, nonmoment-based measures, quantile functions, skewness index.
\color{black}
\section{Introduction}
Asymmetry is a pervasive characteristic in natural phenomena and human decision-making under uncertainty, and its proper quantification is crucial across diverse fields. 
The concept is fundamental, and its effective definitions and applications vary significantly across probabilistic frameworks and applications, but surprisingly, a similar diversity of adaptive definitions, allowing for tailored quantification across diverse scenarios, has not been reflected in fuzzy set theory.

Our focus here is on the need for a robust and interpretable measure of asymmetry of fuzzy numbers within financial engineering applications, particularly in portfolio optimization, where investor preferences for gains versus losses necessitate a clear understanding of skewness.

Modern portfolio theory lends itself as a prime field of application which necessitates a notion of asymmetry to distinguish between the preference for gains and losses. Markowitz~\cite{Markowitz1952} introduced seminal groundwork for this theory in 1952, aiming to formulate an efficient portfolio with a maximal ratio of mean to standard deviation, the Sharpe ratio. 
This framework relies solely on means and variance as representations of portfolio risk and return. 
However, recent advances in decision theory have discovered the importance of higher order moment risk preferences in human decision making~\cite{Menzesetal80,eeckhoudtachlesinger06,trautmanndekuilen18}, oftentimes characterized through the preference for a positive sign of higher order odd moments~\cite{ebert13}. Accordingly, adaptations and extensions of the traditional mean-variance theory have also been proposed in portfolio theory and financial economics \cite{whitmore75, krauslitzenberger76,harveysiddique00}.

These advancements in probabilistic finance \cite{andersenetal03, neuberger13, kozhanneubergerschneider13, bondarenko14b, schneidertrojani18} underscore the financial community's significant interest in robust definitions of higher order return moments, such as skewness and kurtosis, highlighting the essential role of traditional probability theory in this context.  However, situations often arise where the underlying probability distributions are unknown, vague, or imprecise, necessitating frameworks that can handle such epistemic uncertainty differently. While nonadditive probabilities (see \cite{gilboaschmeidler89, Gilboa2009}), ambiguity and Knightian uncertainty in financial economics \cite{epsteinschneider08,anthropelosschneider24} offer theoretical alternatives, their specification is often subject to subjective modeling choices and potential nonrobustness to perturbations of inputs.
This underscores the need for similarly rigorous and interpretable measures within fuzzy contexts when addressing epistemic uncertainty.

Fuzzy sets are highly developed and suitable devices for robustly modeling and analyzing vague or imprecise quantities, as their membership functions can be systematically constructed from observed data or mathematical models~\cite{ma2022data}.
Numerous works highlight the application of fuzzy set theory in financial economics and portfolio theory. This literature is broadly categorized by its foundational theories, primarily \emph{possibilistic}, \emph{credibilistic}, and \emph{intuitionistic} (see Section~\ref{subsec:FuzzyMeanValues}). Early work by Dijkman et al.~\cite{zadeh1983fuzzy} was the starting point for scholarly focus on portfolio selection within a fuzzy framework. Subsequent contributions explore  comparisons between fuzzy and stochastic programming \cite{inuiguchi2000possibilistic}, proposals for possibilistic variances and covariances \cite{zhang2003possibilistic}, analyses of semivariance for fuzzy variables \cite{huang2008mean}, and mean-absolute deviation models using fuzzy linear regression \cite{inuiguchi1997mean}. These efforts underscore the flexibility of fuzzy set theory in handling financial uncertainty.

Despite this rich body of work, the existing literature on fuzzy skewness exclusively relies on direct adaptations of moment-based skewness coefficients from classical probability theory, particularly those akin to the Fisher--Pearson coefficient~\cite{Pearson1895}. These approaches habitually make an implicit, albeit interpretationally unfounded, identification between a fuzzy membership function $\mu(x)$ and a probability density function (PDF) $f_X(x)$ of an underlying random variable $X$, often without explicitly defining the random variable itself.
A membership function's value indicates a degree of belonging or possibility without a direct probabilistic interpretation in a general context. Consequently, simply substituting membership functions into formulas designed for probability density functions leads to definitions lacking rigorous theoretical justification. This philosophical gap often overlooks the fundamental differences in the underlying theories of possibility and probability~\cite{Zad99,Dub06,DB12,HH19,KZ21}. This underscores the need for fuzzy skewness coefficients that are not merely formulaic analogues, but are rigorously derived from the intrinsic properties of fuzzy numbers and their relationship to probabilistic concepts.

This article argues that a more rigorous approach is needed, particularly one that leverages the inherent structure of fuzzy numbers (their $\alpha$-cuts) and establishes a clear, explicit link to probability theory, rather than relying on analogical substitutions.

Thus, we propose a quantile-based skewness measure for fuzzy numbers built upon an explicit probabilistic interpretation. Rather than an entirely new notion of skewness, the proposal adapts well-established quantile-based coefficients from classical statistics, those of Bowley~\cite{bowley1926elements}, Kelley~\cite{kelley1921new}, and Groeneveld and Meeden~\cite{groeneveld1984measuring}, to the $\alpha$-cut representation of fuzzy numbers, and complements, rather than replaces, the existing moment-based definitions. Given the fat tails commonly observed in particular in financial data, where higher order moments may not exist, such a quantile-based definition is desirable for robustness reasons.

The new skewness coefficient is not only symmetric and invariant under transformations of scale and location, but also requires fewer arithmetic operations per evaluation than moment-based skewness coefficients. This approach connects possibility and probability theory by grounding the construction of fuzzy numbers in explicit probabilistic terms, interpreting the left and right components of membership functions as cumulative and survival probability functions, respectively. We thereby provide an interpretation of the $\alpha$-cuts as generalized quantiles, forming the basis for our quantile-based skewness measure.

To illustrate the practical utility of our approach, we operationalize the proposed quantile-based fuzzy skewness measures within a mean–variance–skewness portfolio optimization framework using historical market data. Specifically, we reconstruct fuzzy asset returns directly from the empirical quantile functions of nine highly volatile, large-cap U.S. growth equities over the period 2012–2019. This setting allows us to assess the out-of-sample performance of our model against traditional moment-based benchmarks under diverse market regimes. To support reproducible research, the complete implementation of both proposed indices alongside the optimization frameworks is made available in a companion repository (link withheld during double-anonymous review).

The rest of this article is organized as follows. Section~\ref{sec:preliminaries} revisits mathematical preliminaries of fuzzy numbers; definition, terminology, arithmetic operations, and gives an overview of existing definitions of mean, variance, and skewness. Section~\ref{sec:our-proposal} introduces our new skewness coefficients for fuzzy numbers.
Section~\ref{sec:empirical} discusses the mean-variance-skewness generalization of Markowitz's Modern Portfolio  Theory, and compares our coefficients against others in terms of out-of-sample portfolio performance and analytical computational complexity. Finally, Section~\ref{sec:conclusions} concludes this article and provides future research outlook.

\section{Mathematical Preliminaries}\label{sec:preliminaries}

This section briefly recalls the basic definitions and properties related to fuzzy numbers, alongside prevalent definitions of their mean, variance, and skewness. 

\smallskip
For this article, we adopt the following definition of a fuzzy number, based on~\cite{zadeh1965fuzzy}:

\begin{definition}\label{def:GeneralFN}
A fuzzy number $\xi$ is defined by an ordered pair of continuous functions $\xi_l(\cdot)$ and $\xi_r(\cdot)$ with compact domains $\mathcal{D}(\xi_l)$ and $\mathcal{D}(\xi_r)$ in $\mathbb{R}$ such that the following holds.

\smallskip
\begin{itemize}
    \item $\max\mathcal{D}(\xi_l) \leq \min\mathcal{D}(\xi_r)$.
    \item $\xi_l(x)$ is strictly increasing onto $[0,1]$.
    \item $\xi_r(x)$ is strictly decreasing onto $[0,1]$.
\end{itemize}

\smallskip
\noindent We denote the respective domains of $\xi$ as follows:
\begin{itemize}
    \item $\mathcal{D}(\xi_l) = [l_{\xi}, \underline{m}_{\xi}]$ 
    \item $\mathcal{D}(\xi_r) = [\overline{m}_{\xi}, r_{\xi}]$
\end{itemize}

\smallskip
The functions $\xi_l(\cdot)$ and $\xi_r(\cdot)$ are referred to as the \textit{left} and \textit{right components} of $\xi$, respectively. The interval $[\underline{m}_{\xi}, \overline{m}_{\xi}]$ is called the \textit{core} of the fuzzy number $\xi$.
\end{definition}

\smallskip
\paragraph{Membership Function}
\noindent The component pair $(\xi_l, \xi_r)$ induces a unique continuous \textit{membership function} $\mu_{\xi}(x): \mathbb{R} \to [0,1]$, defined as
\begin{equation}
\mu_{\xi}(x) = \begin{cases}
\xi_l(x) & \text{if } l_{\xi} \leq x < \underline{m}_{\xi},\\
1 & \text{if } \underline{m}_{\xi} \leq x \leq \overline{m}_{\xi},\\
\xi_r(x) & \text{if } \overline{m}_{\xi} < x \leq r_{\xi},\\
0 & \text{otherwise.}
\end{cases}
\end{equation}

\smallskip
\paragraph{Level-set Representation}
\noindent The components $\xi_l(\cdot)$ and $\xi_r(\cdot)$ are invertible by design. For each $\alpha \in [0,1]$, we may denote:

\smallskip
\begin{itemize}
    \item $\xi_d(\alpha) := \xi_l^{-1}(\alpha)$,
    \item $\xi_u(\alpha) := \xi_r^{-1}(\alpha)$.
\end{itemize}

\smallskip
This leads to an equivalent representation of the fuzzy number $\xi$ as an indexed family of closed intervals,
\begin{equation}
C_{\alpha}(\xi) := [\xi_d(\alpha), \xi_u(\alpha)], \quad \alpha \in [0,1].
\end{equation}
We refer to the interval $C_{\alpha}(\xi)$ as the \textit{$\alpha$-cut} of the fuzzy number $\xi$, and to the functions $\xi_d(\alpha)$ and $\xi_u(\alpha)$ as the \textit{lower} and \textit{upper bound functions} of $\xi$, respectively.

\smallskip
\noindent Note that $C_1(\xi) = [\underline{m}_{\xi}, \overline{m}_{\xi}]$ (the core of $\xi$), and $C_0(\xi) = [l_{\xi}, r_{\xi}]$, the \textit{support} of its membership function.

\subsection{Arithmetic Operations}
The level-set representation allows for a convenient definition of arithmetic operations on fuzzy numbers. For two fuzzy numbers $\xi$ and $\eta$, and an arithmetic operation $\circ \in \{\oplus, \ominus, \odot, \oslash\}$, the $\alpha$-cuts of $\xi \circ \eta$ are given by
\begin{equation}
\xi\circ\eta := [{(\xi\circ \eta)}_d, {(\xi\circ \eta)}_u]
\end{equation}
where
\begin{equation}\label{FuzzyArithmetic-left}
{(\xi\circ \eta)}_d:= \min\bigl\{\xi_d\diamond \eta_d, \xi_d\diamond \eta_u, \xi_u\diamond \eta_d, \xi_u\diamond \eta_u\bigr\},
\end{equation}
and
\begin{equation}\label{FuzzyArithmetic-right}
{(\xi\circ \eta)}_u:= \max\bigl\{\xi_d\diamond \eta_d, \xi_d\diamond \eta_u, \xi_u\diamond \eta_d, \xi_u\diamond \eta_u\bigr\},
\end{equation}

\noindent where $\diamond$ denotes the corresponding crisp (of real numbers) arithmetic operation ($+,-,\cdot,/\,;$ for division, $0\not\in supp[\eta]$ is required). 

\smallskip

\subsection{Probabilistic Interpretation of Components}\label{subsec:ProbInterpretation}
As continuous and strictly monotonic functions onto $[0,1]$, the left and right components of a fuzzy number $\xi$ (see Definition~\ref{def:GeneralFN}) satisfy the definitional requirements of a cumulative distribution function and a survival function, respectively.

\smallskip
We may thus associate the left component $\xi_l$ with a continuous random variable $X_L$, whose cumulative distribution function (CDF), $F_{X_{L}}(x) = \Pr(X_{L} \leq x)$, is given as
\begin{equation}\label{eq:leftComponent2leftRandom}
F_{X_{L}}(x) = \begin{cases}
0 & \text{if } x < l_{\xi}, \\
\xi_l(x) & \text{if } l_{\xi} \leq x \leq \underline{m}_{\xi}, \\
1 & \text{if } x > \underline{m}_{\xi}.
\end{cases}
\end{equation}

Similarly, we may associate the right component $\xi_r$ with a continuous random variable $X_R$, given by its survival function $S_{X_{R}}(x) = \Pr(X_R > x)$ as
\begin{equation}\label{eq:rightComponent2rightRandom}
S_{X_{R}}(x) = \begin{cases}
1 & \text{if } x < \overline{m}_{\xi}, \\
\xi_r(x) & \text{if } \overline{m}_{\xi} \leq x \leq r_{\xi}, \\
0 & \text{if } x > r{\xi}.
\end{cases}
\end{equation}

Equivalently, the CDF of of the associated random variable $X_R$ is given by $F_{X_R}(x) = 1 - S_{X_R}(x)$ for $x \in \mathbb{R}$.

\smallskip

\noindent Under this construction, the lower bound $\xi_d(\alpha)$ is the point for which $\Pr(X_L \leq \xi_d(\alpha)) = \alpha$, implying that $X_L$ exceeds this value with probability $1-\alpha$. Similarly, the upper bound $\xi_u(\alpha)$ satisfies $\Pr(X_R > \xi_u(\alpha)) = \alpha$. Consequently, the parametric bound functions of the fuzzy number $\xi$ can be expressed in terms of the \textit{quantile functions} $Q_{X}(\cdot) = F^{-1}_{X}(\cdot)$ of the associated random variables,
\begin{equation}\label{eq:Quantile2Parametric}
\begin{aligned}
& \xi_d(\alpha) = Q_{X_L}(\alpha), \\
& \xi_u(\alpha) = Q_{X_R}(1-\alpha),
\end{aligned}
\end{equation}
which allows the $\alpha$-cuts of $\xi$ to be written as:
\begin{equation}
C_{\alpha}(\xi) = [Q_{X_L}(\alpha), Q_{X_R}(1-\alpha)].
\end{equation}

\smallskip
This gives an interpretation of a fuzzy number's $\alpha$-cuts as representing intervals of values which are ``at least $(1-\alpha)$-probable'' in terms of their respective associated distribution and survival functions, establishing a structured semantic bridge between degrees of possibility (fuzzy membership) and probabilistic quantiles.

\subsection{Review of Fuzzy Mean, Variance, and Skewness Coefficients in the literature}\label{subsec:FuzzyMeanValues}

Building upon the conceptual groundwork laid by Schneider \cite{schneider2020new}, who proposed a geometrically motivated approach to fuzzy skewness and identified the potential for constructing fuzzy numbers from CDFs, the present study formalizes this direction by deriving a comprehensive quantile-based skewness framework and operationalizing it for practical portfolio optimization.

While the mathematical Definition~\ref{def:GeneralFN}, is universally accepted subject to minor variants regarding continuity and differentiability conditions and support, there is a significant diversity in definitions for fuzzy mean, and consequently variance and skewness.

Broadly speaking, these definitions stem from distinct theoretical frameworks interpreting the semantic value of fuzzy numbers. The most prominent frameworks are \emph{possibility theory}~\cite{Zad78,denoeux2023reasoning,yager2012introduction,ma2022data}, \emph{credibility theory}~\cite{LL02,das2023quantum,zhang2018portfolio,das2023pentagonal}, and \emph{intuitionistic fuzzy sets}~\cite{Ata86,Ata89}.

Each framework proposes its respective definitions of mean values. 

The original \emph{possibilistic} mean value is defined in~\cite{GV86} as
\begin{equation}\label{eq:def:PossMean}
\E_{\mathcal{P}}(\xi) := \dfrac{1}{2}\cdot \left(\int_{0}^{1}\xi_d(\alpha)\,d\alpha + \int_{0}^{1}\xi_u(\alpha)\,d\alpha\right).
\end{equation}

\noindent Numerous modifications to this definition exist, often incorporating weighting by $\alpha$~\cite{FM03}. 

\smallskip
The \emph{credibilistic} mean value is given by
\begin{equation}\label{eq:def:CredMean}
\E_{Cr}(\xi) := \int_{l_{\xi}}^{r_{\xi}}x\phi_{Cr}(x)dx,
\end{equation}
where $\phi_{Cr}(x)$ is the credibilistic density function defined as the derivative of the credibilistic distribution of $\xi$ defined by
$$\Phi(x) = \dfrac{1}{2}\left(\sup_{y\leq x}\xi(y) + 1 -\sup_{y>x}\xi(y)\right).$$
See~\cite{ZYW16} for some examples of credibilistic distributions and mean values.

\smallskip
For intuitionistic definitions see~\cite{Ata86, Ata89}.

\smallskip
\begin{remark}
It is easy to see by change of variables that for invertible left and right components the two integrals~\eqref{eq:def:PossMean} and~\eqref{eq:def:CredMean} are the same, i.e.
\begin{equation}
\E_{Cr}(\xi) = \E_{\mathcal{P}}(\xi)
\end{equation}
\end{remark}

Beyond these definitions, the literature offers various fuzzy number \emph{ranking functions} for measuring central tendency (see~\cite{BM13} for an overview).

In parallel and in consequence to the above definitions of central tendency, credibilistic, possibilistic, and intuitionistic variances and skewness coefficients have been developed. The definitions then habitually draw direct analogy from the definition of Fisher--Pearson central moments in probability theory~\cite{Pearson1895}.

Thus, once the definition for the  mean value $\E(\xi)$ has been chosen an obvious possible definition for fuzzy variance may consequently be
\begin{equation}\label{eq:FuzzyVarianceCleanPearson}
\V(\xi) = \E(\xi - \E(\xi))^2,
\end{equation}
where the membership function $\xi(x)$ replaces in a straightforward manner the probability density in the probabilistic definition.

Analogously, fuzzy skewness can be defined in analogy to the probabilistic Fisher--Pearson coefficient along the lines of
\begin{equation}\label{eq:FuzzySkewnessCleanPearson}
\Skew(\xi) = \E(\xi - \E(\xi))^3/\V(\xi)^{\frac{3}{2}},
\end{equation}
where the expected value $\E$ in the above definitions is chosen to be credibilistic, possibilistic, and intuitionistic, or other.

Importantly, such literal translation preserves the invariance property~\cite{Zwe64}: for $c>0$ and $d \in \mathbb{R}$, $\Skew(c\,\xi+d) = \Skew(\xi)$, meaning that skewness so defined is invariant under transformations of scale and location, provided that the chosen mean and variance operators are themselves appropriately invariant.

However, authors chose to alter /modify the definitions~\eqref{eq:FuzzyVarianceCleanPearson} and~\eqref{eq:FuzzySkewnessCleanPearson} to better suit intuition and/or interpretational value within the framework of fuzzy set theory, and thereby at the same time potentially compromise mathematical consistency.
For example: In the representative examples that we have chosen to compare our new skewness coefficient against, 
Li et al.~\cite{li2015skewness} adhered to $\Skew(\xi) = \E(\xi - \E(\xi))^3$  (without any division by a power of variance), while Vercher--Bermudez~\cite{VB13} are of type $\E(\xi - \E(\xi))^3/\V(\xi)^{\frac{3}{2}},$ albeit with the Fisher--Pearson type standard deviation replaced by the dimensionally \emph{but not asymptotically} equivalent measure \emph{downside risk} first introduced in~\cite{VB09}. 

The extensive literature on fuzzy portfolio optimization involving a skewness constraint further exemplifies this practice, with contributions categorized by their theoretical tradition: possibilistic \cite{VB09, LZX12, VB13, HB15, YWD15, SRBVL16, LZZ18, LY19, MH20}, credibilistic \cite{CL10, KDF12, BAM13, XLLL17, BKKM09, PJ21, LWLWH21, MGK21}, and intuitionistic \cite{CLLYY11, DP18, GMYK19, YKMGC23}. 

Still, all fuzzy approaches mentioned rely on probabilistic moment-based skewness, directly substituting fuzzy membership functions for probability density functions in probabilistic formulae. None consider classical alternative quantile-based skewness coefficients (e.g., Bowley’s~\cite{bowley1926elements}, Kelley's~\cite{kelley1921new}, Pearson’s mode or median skewness) or modern nonmoment-based approaches of contemporary mathematical finance and financial engineering~\cite{KL73, Neu12, CET13, Ebe13, BKM03, ST19,  SWZ20, NP21}.

\subsection{Skewness Coefficients With Which We Compare}\label{subsec:ComparisonCoefficients}
We select two prominent fuzzy skewness coefficients from the literature for comparison: those proposed by Vercher--Bermudez~\cite{VB13} in 2013 and by Li et al.~\cite{li2015skewness} in 2015.

\subsubsection{Vercher--Bermúdez Coefficient (2013)}\label{subsubsec:VB13_skewness}
Vercher and Bermúdez~\cite{VB13} define mean, variance (actually \emph{downside risk}, $\omega(\xi)$, which is dimensionally equivalent to standard deviation), and skewness of fuzzy numbers as follows,

\begin{align}
	\E_{VB13}(\xi) &= 1/2 \cdot \int_{0}^{1}[\xi_d(\alpha) + \xi_u(\alpha)]\,d\alpha \quad \text{(as in \eqref{eq:def:PossMean})},\label{eq:BermMean}\\
	(\V_{VB13})^{\frac{1}{2}} &= \omega(\xi) := \int_{0}^{1}[\xi_u(\alpha) - \xi_d(\alpha)]\,d\alpha.\label{eq:VarBermudez}
\end{align}
The \emph{third possibilistic moment} is defined as
\begin{align}\label{eq:BermMu3}
	\mu_3(\xi) &= 1/2\cdot\int_{0}^{1}[\xi_d(\alpha) - \E(\xi)]^3\,d\alpha \\
	&\quad + 1/2\cdot\int_{0}^{1}[\xi_u(\alpha) - \E(\xi)]^3\,d\alpha . \notag
\end{align}

The final formula for skewness is ratio of the third possibilistic moment and the third power of downside risk,

\begin{equation}\label{eq:BermSkew}
	\Skew_{VB13}(\xi) = \dfrac{\mu_3(\xi)}{[\omega(\xi)]^3}.
\end{equation}

\subsubsection{Li--Guo--Yu Coefficient (2015)}\label{subsubsec:LGY15_skewness}
The definitions for mean, variance, and skewness introduced by Li et al.~\cite{li2015skewness} are given by the following formulae,

\begin{align}
	\E_{LGY15}(\xi) &= \int_{0}^{1}\alpha\cdot[\xi_d(\alpha) + \xi_u(\alpha)]\,d\alpha , \\\label{eq:LiMean}
	\V_{LGY15}(\xi)
		&= \int_{0}^{1}\alpha \cdot \Bigl[(\xi_d(\alpha) - \E(\xi))^2 \\ \notag
		&\quad\qquad\qquad + (\xi_u(\alpha) - \E(\xi))^2\Bigr]\,d\alpha,  \\ \label{eq:LiVar}
	\Skew_{LGY15}(\xi)
		&= \int_{0}^{1}\alpha \cdot \Bigl[(\xi_d(\alpha) - \E(\xi))^3 \\ \notag
		&\quad\qquad\qquad + (\xi_u(\alpha) - \E(\xi))^3\Bigr]\,d\alpha. \label{eq:LiSkew}
\end{align}

\section{Our Proposal}\label{sec:our-proposal}
Given the plethora of moment-based definitions for fuzzy mean, variance and skewness the nonmoment-based (that is not explicitly incorporating mean and variance) quantile-based skewness coefficients proposed by Bowley~\cite{bowley1926elements}, Kelley~\cite{kelley1921new}, and in particular by Groeneveld and Meeden~\cite{groeneveld1984measuring} offer a viable alternative for adaption to fuzzy numbers.

The skewness coefficient of Groeneveld and Meeden, who relate to the theoretical foundations presented by van Zwet in~\cite{Zwe64}, 
has instigated a multitude of further theoretical research and has also lent itself to various practical applications e.g.~\cite{langlois2020measuring,dai2021skewness}.

\smallskip
The Groeneveld and Meeden~\cite{groeneveld1984measuring} skewness coefficient is defined for a random variable $X$ as:

\begin{equation}\label{eq:GM_Skewness1}
\Skew_{GM_1}(X)(p) := \frac{Q_X(1-p) + Q_X(p)- 2\,Q_X(1/2)}{Q_X(1-p) - Q_X(p)},
\end{equation}
for $p \in \left(0, 1/2\right)$, where $Q_X(\cdot)$ denotes the quantile function of the random variable $X$.

Bowley's coefficient~\cite{bowley1926elements} is a special case of~\eqref{eq:GM_Skewness1} and occurs at $p = 0.25$, and Kelley's coefficient~\cite{kelley1921new} at $p = 0.1.$

A remark on terminology: in classical statistics the \emph{quartiles} are the particular quantiles at levels $p = 0.25, 0.5, 0.75$, whereas the \emph{quantile function} $Q_X(p)$ is defined for arbitrary $p \in (0,1)$. Our coefficients are defined for arbitrary levels $p$ (equivalently, arbitrary $\alpha$-levels of a fuzzy number), with quartile-based measures such as Bowley's arising only as the special case $p = 0.25$; hence we call the proposed measures --- and this article's title --- \emph{quantile}-based rather than \emph{quartile}-based.

In the same paper~\cite{groeneveld1984measuring}, Groeneveld and Meeden~\cite{groeneveld1984measuring} proposed a second skewness coefficient by taking the integral of numerator and denominator of the first~\eqref{eq:GM_Skewness1},

\begin{equation}\label{eq:GM_Skewness2}
\Skew_{GM_2}(X) := \frac{\int_{0}^{0.5} (Q_X(1-p) + Q_X(p) - 2 Q_X(0.5))\, dp}{\int_{0}^{0.5} (Q_X(1-p) - Q_X(p))\, dp}.
\end{equation}

\subsection{Definitions of the New Skewness Coefficients}

Let $\xi$ be a fuzzy number with core $[\underline{m}_{\xi}, \overline{m}_{\xi}]$ and $\alpha$-cuts $[\xi_d(\alpha), \xi_u(\alpha)]$. Let $X_L$ and $X_R$ be the random variables associated with the left and right components of $\xi$, as per~\eqref{eq:leftComponent2leftRandom} and \eqref{eq:rightComponent2rightRandom}. Based on \eqref{eq:Quantile2Parametric}, their quantile functions are given by $Q_{X_L}(\alpha) = \xi_d(\alpha)$ and $Q_{X_R}(\alpha) = \xi_u(1-\alpha)$, respectively.

In analogy to the two skewness coefficients of Groeneveld and Meeden, the point-based $\Skew_{GM_1}$ and the integrated $\Skew_{GM_2}$ we also define two coefficients, the first one point-based, and the second based on the first one by integration.

\smallskip
\begin{definition}\label{def:OurFirstCoefficient}$ $\par\noindent
The point-based fuzzy skewness coefficient $\Skew_{{JKPT}_1}(\xi)(\alpha)(v)$ is defined as the convex combination of two constituents, pertaining to what we choose to call the \emph{inner skewness}, $$\Skew_{Inner,point}(\xi)(\alpha),$$
and the \emph{outer} skewness, $$\Skew_{Outer,point}(\xi)(\alpha),$$
of a fuzzy number $\xi$ at point~$\alpha,$ i.e.
\begin{multline}\label{eq:OurFirstCoefficient}
\Skew_{{JKPT}_1}(\xi)(\alpha)(v) := 
v \cdot \Skew_{Outer,point}(\xi)(\alpha)\\ 
+ (1-v) \cdot \Skew_{Inner,point}(\xi)(\alpha), \, v\in[0,1].
\end{multline}
Here, the following holds.
\smallskip
\begin{enumerate}
    \item \textbf{Point-based Inner Skewness ($\Skew_{Inner,point}(\xi)(\alpha)$):}
    This constituent captures the intrinsic skewness of the underlying probabilistic distributions $X_L$ and $X_R$. For a chosen $\alpha$-level $\alpha \in (0, 0.5)$, it is defined as the average of the Groeneveld - Meeden skewness coefficients for $X_L$ and $X_R$, expressed by identifications \eqref{eq:leftComponent2leftRandom} and \eqref{eq:rightComponent2rightRandom} in terms of the $\alpha$-cut functions $\xi_d$ and $\xi_u.$ 
    
    \begin{equation}\label{eq:S_InnerGM1}
    \Skew_{inner,point}(\xi) (\alpha) := 
    \frac{1}{2} \left( S_{GM_1}(\xi_d)(\alpha) + S_{GM_1}(\xi_u)(\alpha) \right),
    \end{equation}
    where 
    \begin{equation}\label{eq:S_InnerLeft1}
    \Skew_{GM_1}(\xi_d)(\alpha) := \frac{\xi_d(1-\alpha) + \xi_d(\alpha)- 2\,\xi_d(1/2)}{\xi_d(1-\alpha) - \xi_d(\alpha)}
    \end{equation}
    for $X_L$, and for $X_R$,
    \begin{equation}\label{eq:S_Right1}
    \Skew_{GM_1}(\xi_u)(\alpha) := \frac{\xi_u(\alpha) + \xi_u(1-\alpha)- 2\,\xi_u(1/2)}{\xi_u(\alpha) - \xi_u(1-\alpha)}.
    \end{equation}
    
    A typical choice for $\alpha$ will be analogous to the probabilistic setting $\alpha = 0.25$ (Bowley's coefficient) or $\alpha= 0.1$ (Kelley's coefficient). 
    
    \bigskip
    \item \textbf{Point-based Outer Skewness ($\Skew_{Outer,point}(\xi)(\alpha)$):}
    This constituent captures the perceived skewness of the fuzzy number's overall profile using a specific $\alpha$-level $\alpha \in (0, 1)$. It is defined as
    \begin{equation}\label{eq:S_Outer_point}
    \Skew_{Outer,point}(\xi)(\alpha) := \frac{(\xi_d(\alpha) - \underline{m}_{\xi}) + (\xi_u(\alpha) - \overline{m}_{\xi} )}{(\xi_u(\alpha) - \xi_d(\alpha))}.
    \end{equation}
\end{enumerate}
\end{definition}

The coefficient's two parts are reminiscent of the law of total variance: $\Skew_{Inner,point}(\xi)(\alpha)$ accounts for the intrinsic skewness of the fuzzy number's left and right components, reflecting the inherent asymmetry of its probabilistic 'building blocks.' $\Skew_{Outer,point}(\xi)(\alpha)$ captures the ``outer'' skewness, reflecting the overall shape's asymmetry via $\alpha$-cuts. This part is particularly relevant for visual assessment and fuzzy numbers whose components might be symmetric but whose overall shape appears skewed (e.g., triangular fuzzy numbers \cite{SKM23, KGMS21}).
This coefficient is computationally very efficient as it only requires evaluating the $\alpha$-cut functions at a few specific points.

\smallskip
In analogy to~\eqref{eq:GM_Skewness2} we introduce a second fuzzy skewness coefficient, $\Skew_{JKPT_2}(\xi)$, which utilizes integral forms for both its inner and outer constituents.

\begin{definition}\label{def:OurSecondCoefficient}$ $\par\noindent
The second  integral-based fuzzy skewness coefficient -
$\Skew_{{JKPT}_2}(\xi)$ is defined as
\begin{equation}\label{eq:OurSecondCoefficient}
\begin{split}
\mathbf{\Skew_{{JKPT}_2}}(\xi)(v) &:= v \cdot \Skew_{Outer,integral}(\xi)\\ 
&+ (1-v) \cdot \Skew_{Inner,integral}(\xi), v \in[0,1].
\end{split}
\end{equation}
Here, the following holds.
\smallskip 
\begin{enumerate}
    \item \textbf{Integral inner skewness ($\Skew_{Inner,integral}(\xi)$):}
    This constituent integrates the intrinsic skewness of the underlying probabilistic distributions $X_L$ and $X_R$ across $\alpha$-levels from $0$ to $0.5$. It is defined as the average of the integral Groeneveld--Meeden skewness coefficients for $X_L$ and $X_R$, expressed in terms of $\alpha$-cut functions,
    \begin{equation}\label{eq:S_Inner_Integral}
    \Skew_{Inner,integral}(\xi) := \frac{1}{2} \left( S_{GM_2}(\xi_d) + S_{GM_2}(\xi_u) \right),
    \end{equation}
    where 
    \begin{equation}
    \Skew_{GM_2}(\xi_d) := \frac{\int_{0}^{0.5} \left(\xi_d(1-\alpha) + \xi_d(\alpha) - 2 \xi_d(0.5)\right) d\alpha}{\int_{0}^{0.5} (\xi_d(1-\alpha) - \xi_d(\alpha)) d\alpha},
    \end{equation} for $X_L$, and for $X_R$,
    \begin{equation}
    \Skew_{GM_2}(\xi_u) := \frac{\int_{0}^{0.5}\, (\xi_u(\alpha) + \xi_u(1-\alpha) - 2 \xi_u(0.5))\, d\alpha}{\int_{0}^{0.5}\, (\xi_u(\alpha) - \xi_u(1-\alpha))\, d\alpha}.
    \end{equation}
    \item \textbf{Integral outer skewness ($\Skew_{Outer,integral}(\xi)$):}
    This constituent captures the emergent skewness of the fuzzy number's overall profile by integrating a quantile-based skewness measure across $\alpha$-levels from $0$ to $0.5$. It is defined as the ratio of integrals,
    \begin{multline}\label{eq:S_Outer_integral}\noindent
    \Skew_{Outer,integral}(\xi) := \\
    \frac{\int_{0}^{0.5}\, (\xi_d(\alpha) - \underline{m}_{\xi}) + (\xi_u(\alpha) - \overline{m}_{\xi})\, d\alpha}{\int_{0}^{0.5}\, (\xi_u(\alpha) - \xi_d(\alpha))\, d\alpha}.
    \end{multline}
\end{enumerate}
\end{definition}

The choice between $\Skew_{JKPT_1}$ and $\Skew_{JKPT_2}$ depends on the desired balance between computational efficiency and comprehensiveness. $\Skew_{JKPT_1}$ (with point-based inner and outer constituents) offers superior speed, while $\Skew_{JKPT_2}$ (with integral inner and outer constituents) provides a more holistic view of the component distributions' skewness. The weighting factor $v \in [0,1]$ allows for customization; for instance, $v=0.5$ provides an equal emphasis on both the emergent fuzzy shape and the intrinsic component distributions. The choice of $\alpha$ for $\Skew_{Inner,point}(\xi) (\alpha)$ and $\alpha$ for $\Skew_{Outer,point}(\xi)(\alpha)$ also influences the coefficient's sensitivity to central versus tail asymmetry (e.g., $\alpha = 0.1$ or $\alpha = 0.25$). These choices should be guided by the specific application context and desired analytical focus.

\smallskip 
\textbf{Notation:} For $v=0.5$ we may omit $v$ in our notation, i.e. $\Skew_{JKPT_1}(\alpha)$ is shorthand for $\Skew_{JKPT_1}(\alpha)(0.5).$ $\Skew_{JKPT_2}$ is shorthand for $\Skew_{JKPT_2}(0.5).$

\subsection{Numerical Examples}\label{subsec:JKPT_Skewness_Examples}
The interpretations outlined in \eqref{eq:leftComponent2leftRandom}--\eqref{eq:Quantile2Parametric} allow for two different starting points to establish the skewness of a given fuzzy quantity: either directly from a fuzzy number which is already given, or conversely, by constructing a fuzzy number from suitably chosen random variables. 

We stress that Definitions~\ref{def:OurFirstCoefficient} and~\ref{def:OurSecondCoefficient} require only the level-set representation of Definition~\ref{def:GeneralFN}. They are not restricted to triangular or otherwise parametric membership functions; accordingly, the examples below employ nontriangular fuzzy numbers (piecewise linear and Beta-generated).
Below we provide two illustrative paragraphs \ref{subsubsec:FN_to_RV} and \ref{subsubsec:RV_to_FN} comprising two examples.

The first Example \ref{ex:FN_to_RV} demonstrates the application of our skewness coefficients - $S_{JKPT_1}(0.1)$, $S_{JKPT_1}(0.25)$, and $S_{JKPT_2}$ as well as $\Skew_{VB13}$ \cite{VB13} and $\Skew_{LGY15}$ \cite{li2015skewness} - to a piecewise linear fuzzy number. This type of fuzzy number typically arises from empirical data (e.g. from histograms) or as a linear approximation of more complex membership functions.

The second example (\ref{ex:RV_to_FN_ExtremeSkew}) takes the reverse approach, beginning with two predefined random variables $X_L$ and $X_R$. From their quantile functions, we construct the left and right components ($\xi_l$ and $\xi_r$) of a fuzzy number, and then apply Definitions~\ref{def:OurFirstCoefficient} and~\ref{def:OurSecondCoefficient}. A step-by-step visualization example of this constructive direction, using moderate magnitudes, is provided in the Supplementary Material for brevity. \par\noindent

\smallskip
\subsubsection{Inferring Generative Random Variables From a Given Fuzzy Number}\label{subsubsec:FN_to_RV}
\begin{example}\label{ex:FN_to_RV}
In this example, a piecewise linear membership function is chosen, a common representation for fuzzy numbers derived from empirical data~\cite{devi1985estimation} or as linear approximations~\cite{klir1995fuzzy},\cite{schneider2018fuzzy} of continuous fuzzy numbers. The fuzzy number $\xi$ is defined by its left component $\xi_l$ and right component $\xi_r$. For reproducibility, we list the transition points $[x, \alpha]$ for each component:
\begin{itemize}
    \item Left component $\xi_l$: 
    $[0.1,0],$ $[0.3,0.1],$ $[0.4,0.2],$ $[0.6,0.25],$ $[0.8,0.4],$ $[1.4,0.5],$ $[1.6,0.7],$ $[2.5,0.85],$ $[2.8,0.9],$ and $[3.0,1]$
    \item Right component $\xi_r$: 
    $[7.6,0],$ $[7.2,0.1],$ $[7.1,0.2],$ $[6.3,0.25],$ $[6.0,0.4],$ $[5.1,0.5],$ $[4.8,0.6],$ $[3.7,0.75],$ $[3.4,0.9],$ and $[3.0,1].$
\end{itemize}

\smallskip
The calculated skewness coefficients for this fuzzy number are: $\Skew_{JKPT_{1}}(0.1) = 0.165,$ $\Skew_{JKPT_{1}}(0.25) = 0.002,$ $\Skew_{JKPT_{2}} = 0.099$ $\Skew_{VB13} = 0.069,$ and $\Skew_{LGY15} = 2.024.$

Fig.~\ref{Fig:FN_to_RV} shows this fuzzy number and its underlying probabilistic interpretation, i.e. CDFs and PDFs of associated random variables whose quantile functions give our skewness coefficient.


\begin{figure}[ht]
  \centering
  \includegraphics[width=0.7\columnwidth]{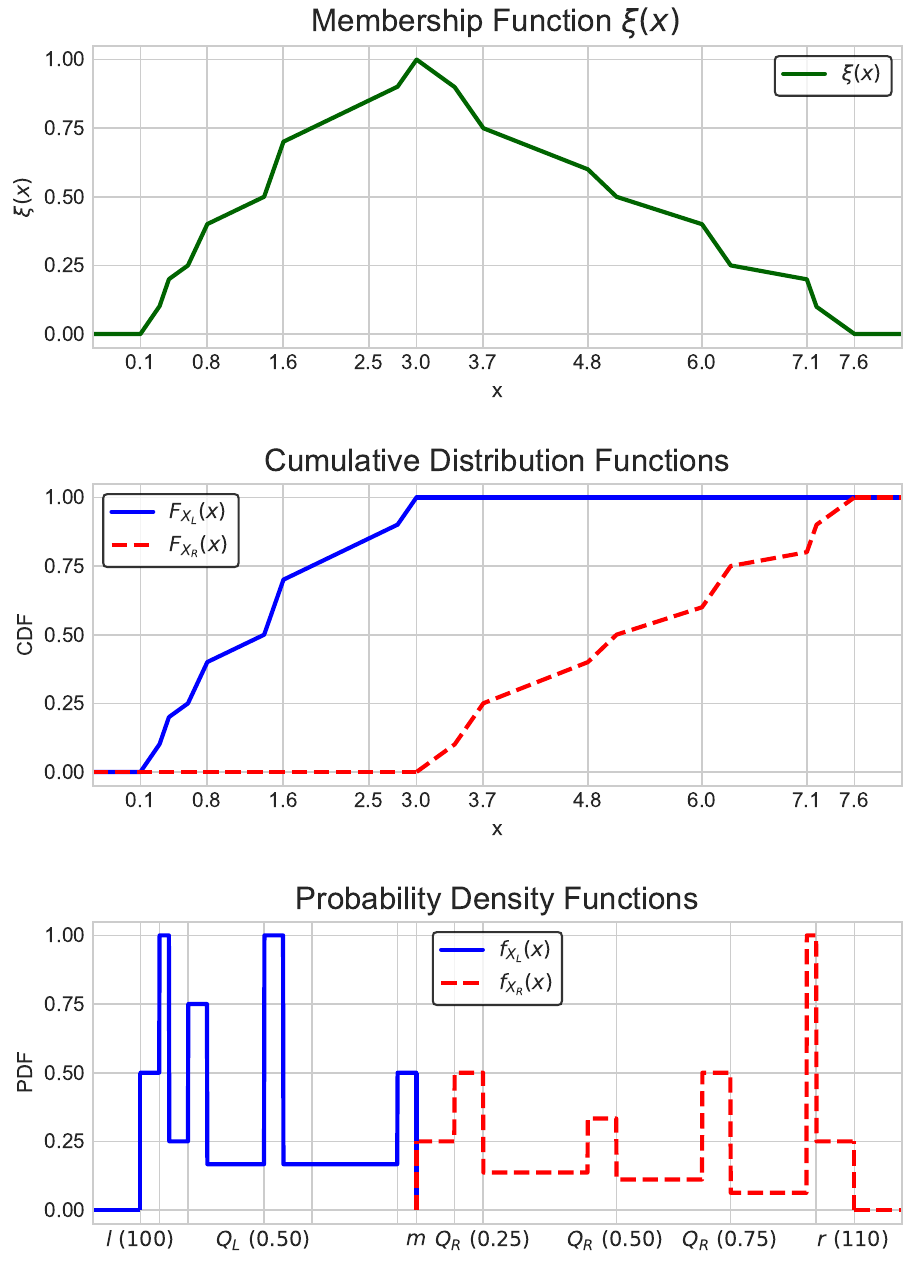}
   \caption{Visualization for Example~\ref{ex:FN_to_RV}:
\newline\textbf{Top:} Piecewise fuzzy number $\xi$.
\newline\textbf{Middle:} Left component $\xi_l$ interpreted as the CDF of $X_L$, and right component $\xi_r$ interpreted as the survival function ($\Pr\{X_R > x\}$) of $X_R$.
\newline\textbf{Bottom:} PDFs obtained by differentiation, with chosen quantile points highlighted.}\label{Fig:FN_to_RV}
\end{figure}
\end{example}

\bigskip
\subsubsection{Constructing a Fuzzy Number From Theoretical Distributions}\label{subsubsec:RV_to_FN}

To illustrate the methodology of construction and properties of the proposed coefficients, we present a stress-test example, designed to showcase the behavior of skewness indices under conditions of extreme asymmetry and larger magnitude, a realistic scenario in financial engineering (e.g., return on Investment for large portfolios) that is often overlooked in fuzzy literature which typically relies on small ``toy'' numbers (e.g. supports like $[1, 10]$). The example also highlights the deficiencies of existing measures like $\Skew_{LGY15}$ and $\Skew_{VB13}$ regarding scale invariance and boundedness. A complementary example visualizing the constructive procedure step by step, transforming two Beta-distributed random variables of moderate magnitude into a single fuzzy number, has been moved to the Supplementary Material for brevity.

\begin{example}{\textit{Extreme Skewness and Scale Invariance}.}\label{ex:RV_to_FN_ExtremeSkew}
A realistic investment scenario often involves large base values with potential for significant, albeit unlikely, upside. Consider a budget or asset value with a base of $100{,}000$ units. We model a scenario where the downside is tight (limited to $1$ unit), but the upside "softness" or potential return stretches significantly ($10\%$ of the base, i.e., $10{,}000$ units). We define $(X_L, X_R)$ as $X_L \sim Beta(0.5, 2)$ on $[100{,}000, \, 100{,}001]$ and $X_R \sim Beta(0.1, 2)$ on $[100{,}001, \, 110{,}000],$ resulting in a fuzzy number with a massive rightward stretch compared to its core.

\smallskip
The calculated skewness coefficients of $\xi$ are: $\Skew_{JKPT_{1}}(0.10)(0.5) = 0.888,$  $\Skew_{JKPT_{1}}(0.25)(0.5) = 0.824,$  $\Skew_{JKPT_{2}}(0.5) = 0.870,$ $\Skew_{VB13} = 34.6,$ $\Skew_{LGY15} \approx 1.5\times 10^{8}$.

Our proposed quantile-based measures yield values of $\approx 0.82$ and $0.87$. These are close to the theoretical maximum of $1.0$, correctly reflecting the extreme right-skewed geometry, and they remain bounded and interpretable. \par
The $\Skew_{LGY15}$ result ($\approx 1.5 \times 10^8$) is very difficult to interpret. This explosion occurs because the definition involves cubic distances from the mean ($\E[(x-\mu)^3]$) but lacks the normalizing division by the cube of the standard deviation ($\sigma^3$) found in standard probability theory. It is not scale-invariant; as the support magnitudes grow (say, from order $10^2$ to order $10^5$), the index grows cubically, making it unsuitable for comparing assets of different magnitudes.

\smallskip
This example illustrates that for financial engineering applications where values can be large and supports variable, a scale-invariant, bounded metric such as $\Skew_{JKPT}$ provides decision support that remains comparable across magnitudes. A formal asymptotic analysis of the ratios between the moment-based indices and the proposed ones, under increasing scale and one-sided spread, is provided in the Supplementary Material.
\end{example}

\subsection{Properties of the New Skewness Coefficient}

\noindent For the following properties, we assume that the fuzzy number $\xi$ is nondegenerate and that the chosen $\alpha$-levels are admissible, ensuring that the relevant interquantile and support spreads (e.g., $\xi_u(\alpha) - \xi_d(\alpha)$) are strictly positive so that the coefficients are well-defined.

The proposed skewness coefficients $S_{JKPT_1}$ and $S_{JKPT_2}$ (beneath we use $S_{JKPT}$ in reference to both) possess several desirable mathematical properties that ensure robustness, interpretability, and suitability for fuzzy numbers:

\begin{enumerate}
    \item Range: $S_{JKPT}(\xi)$  lies within the interval $[-1, 1]$. This is a direct consequence of both the inner ($\Skew_{Inner,point}$ or $\Skew_{Inner,integral}$) and outer ($\Skew_{Outer,point}$ or $\Skew_{Outer,integral}$) constituents individually having a theoretical range of $[-1, 1]$, and the weighting factor $v \in [0, 1]$.
    \item Symmetry: If $\xi$ is a perfectly symmetric fuzzy number (i.e., $\underline{m}_{\xi} - \xi_d(\alpha) = \xi_u(\alpha) - \overline{m}_{\xi}$ for all $\alpha$), then both $\Skew_{Inner,point}$ (or $\Skew_{Inner,integral}$) and $\Skew_{Outer,point}$ (or $\Skew_{Outer,integral}$) are zero, resulting in $S_{JKPT}(\xi) = 0$.
    \item Scale Invariance: $S_{JKPT}(\xi)$ is invariant under positive scaling transformations. For $\xi' = k \xi$ with $k > 0$, $S_{JKPT}(\xi') = S_{JKPT}(\xi)$. This also means the measure is independent of the units of measurement.
    \item Location invariance (translation invariance): $S_{JKPT}(\xi)$ is invariant under additive transformations. For $\xi' = \xi + c$ with $c \in \mathbb{R}$, $S_{JKPT}(\xi') = S_{JKPT}(\xi)$. This ensures that shifting the fuzzy number along the x-axis does not alter its skewness.
    \item Monotonicity (direction of skewness):
        \begin{itemize}
            \item $S_{JKPT}(\xi) > 0$ indicates right-skewness (positive skewness).
            \item $S_{JKPT}(\xi) < 0$ indicates left-skewness (negative skewness).
            \item symmetry implies $S_{JKPT}(\xi) = 0$.
        \end{itemize}
    \item Decomposition and interpretability: The coefficient's formulation as a weighted sum of two constituents offers natural interpretability:
        \begin{itemize}
            \item $\Skew_{Inner,point}(\xi)(\alpha)$ (respectively $\Skew_{Inner,integral}(\xi)$) captures the intrinsic skewness of the underlying probabilistic generative processes ($X_L$ and $X_R$). This provides a deeper, constituent-level understanding of why the fuzzy number exhibits its particular shape. For fuzzy numbers whose components imply symmetric distributions our coefficients give $0.$
            \item $\Skew_{Outer,point}(\xi)(\alpha)$ (respectively $\Skew_{Outer,integral}(\xi)$) captures the perceived skewness of the fuzzy number's overall shape, as seen \par\noindent 
                through  its $\alpha$-cut profile. 
                This constituent is crucial for assessing visual asymmetry, especially in cases where the inner skewness might be zero but the fuzzy number visually appears skewed.
            \item The weighting factor \textcolor{green}{$v$} allows for emphasizing either the perceived (outer) or generative (inner) aspects of skewness, providing flexibility in application and analysis.
        \end{itemize}
    \item Computational efficiency: For point-based versions ($\Skew_{Inner,point}$ and $\Skew_{Outer,point}$), the coefficient is computationally very efficient, involving only a few evaluations of the $\alpha$-cut function. This makes it well suited to settings requiring many repeated evaluations; the scope of this efficiency claim is discussed in Section~\ref{subsec:complexity}. Even the integral forms ($\Skew_{Outer,integral}$ and $\Skew_{Inner,integral}$) involve integration over $\alpha$, which is generally well-behaved for fuzzy numbers and efficient for numerical computation.
\end{enumerate}



\section{Empirical Illustration: Quantile-Based Fuzzy Portfolio Optimization}
\label{sec:empirical}

We illustrate the proposed quantile-based skewness measures in a fuzzy mean--variance--skewness portfolio selection problem with real financial data. The goal is to compare portfolios optimized under three different fuzzy skewness definitions and assess their out-of-sample performance in standard probabilistic terms.
We construct fuzzy assets from the empirical quantile functions of the assets' historical monthly returns.

\subsection{Data and Construction of Fuzzy Asset Returns}
\label{subsec:data}

To ensure the presence of pronounced asymmetry and fat tails, essential for rigorously testing skewness-aware models, we deliberately select a universe of nine highly volatile, large-cap U.S. equities, predominantly from the technology, and growth sector (Apple Inc.\ (AAPL), Microsoft Corporation (MSFT), Amazon.com, Inc.\ (AMZN), NVIDIA Corporation (NVDA), Tesla, Inc.\ (TSLA), Netflix, Inc.\ (NFLX), Advanced Micro Devices, Inc.\ (AMD), Alphabet Inc.\ (GOOGL), and Meta Platforms, Inc.\ (META)). Monthly adjusted closing prices are obtained from Yahoo Finance over the period January 2012 to December 2019, and we compute simple monthly returns. This period covers both benign and volatile regimes, including the 2015--2016 selloff, while avoiding the extreme distortions of the COVID-19 crisis.

Two data details deserve note. First, META (formerly Facebook) was listed only in May 2012, so its return series begins in June 2012 and the panel is slightly unbalanced at the start; all remaining series cover the full period. Second, to avoid look-ahead bias, the empirical quantile functions defining the fuzzy asset returns, and the statistics in Tables~\ref{tab:asset_statistics_competitors} and~\ref{tab:asset_statistics_jkpt}, are estimated on the training window only (the first two thirds of the sample, cf.~Section~\ref{subsec:setup}); the remaining third is reserved for the out-of-sample evaluation.

For each asset, we construct a fuzzy number for its return using empirical quantiles. For an $\alpha$-level $\alpha \in [0,1]$, we define the $\alpha$-cut of the fuzzy return $\xi_i$ of asset $i$ as
\begin{equation}\label{eq:FuzzyAssetConstruction}
C_\alpha(\xi_i) 
=
\big[\xi_{i,d}(\alpha),\,\xi_{i,u}(\alpha)\big] 
=
\big[ Q_{R_i}(\alpha/2),\, Q_{R_i}(1-\alpha/2)\big],
\end{equation}
where $Q_{R_i}$ denotes the \emph{empirical} quantile function of the historical monthly returns of asset $i$. Thus, the lower and upper $\alpha$-cut functions $\xi_{i,d}(\alpha)$ and $\xi_{i,u}(\alpha)$ coincide with the lower and upper empirical quantiles of the return distribution. 

This yields a quantile-based level-set representation entirely consistent with the probabilistic interpretation introduced in Section~\ref{sec:our-proposal}: the left and right components of the membership function correspond to cumulative and survival quantiles, and $\alpha$-cuts can be read as generalized ``at least $(1-\alpha)$-probable'' return intervals: by construction \eqref{eq:FuzzyAssetConstruction}, $\Pr\big(R_i \in C_\alpha(\xi_i)\big) = 1-\alpha$ under the empirical return distribution.

\begin{table}[htbp]
\centering
\caption{Asset Statistical Attributes: Mean, Variance, and Competitor Skewness. Values for $\Skew_{LGY15}$ are near zero due to the rigid cubic deviation from the mean, limiting its effectiveness for fuzzy numbers constructed via symmetric empirical quantiles.}
\label{tab:asset_statistics_competitors}
\footnotesize 
\setlength{\tabcolsep}{4pt}
\begin{tabular}{l rrrrrr}
\toprule
Asset & $\E_{VB13}$ & $\V_{VB13}$ & $\E_{LGY15}$ & $\V_{LGY15}$ & $\Skew_{VB13}$ & $\Skew_{LGY15}$ \\ 
\midrule
AAPL  & 0.017 & 0.013 & 0.019 & 0.002 & -0.050 & -0.000017 \\
MSFT  & 0.017 & 0.007 & 0.017 & 0.001 &  0.106 & -0.000005 \\
AMZN  & 0.028 & 0.015 & 0.027 & 0.002 &  0.050 & -0.000008 \\
NVDA  & 0.035 & 0.016 & 0.029 & 0.002 &  0.221 &  0.000066 \\
TSLA  & 0.047 & 0.052 & 0.034 & 0.007 &  0.498 &  0.000233 \\
NFLX  & 0.045 & 0.048 & 0.038 & 0.006 &  0.473 &  0.000052 \\
AMD   & 0.024 & 0.061 & 0.017 & 0.009 &  0.069 &  0.000381 \\
GOOGL & 0.020 & 0.007 & 0.016 & 0.001 &  0.226 &  0.000012 \\
META  & 0.033 & 0.021 & 0.029 & 0.003 &  0.275 &  0.000036 \\
\bottomrule
\end{tabular}
\end{table}

\begin{table}[htbp]
\centering
\caption{Asset Statistical Attributes: $\Skew_{JKPT}$ Inner and Outer Skewness Components. The proposed indices capture the tail asymmetry of the empirical return distributions under the same median-anchored construction.}
\label{tab:asset_statistics_jkpt}
\footnotesize 
\setlength{\tabcolsep}{4pt}
\begin{tabular}{l rrrrrr}
\toprule
Asset & $\Skew^{1}_{in}(0.1)$ & $\Skew^{1}_{out}(0.1)$ & $\Skew^{1}_{in} (0.25)$ & $ \Skew^{1}_{out} (0.25)$ & $\Skew^{2}_{in}$ & $ \Skew^{2}_{out}$ \\ 
\midrule
AAPL  & -0.230 & -0.028 & -0.179 &  0.016 & -0.202 &  0.039 \\
MSFT  &  0.009 & -0.063 & -0.134 & -0.193 &  0.015 & -0.052 \\
AMZN  &  0.119 &  0.052 &  0.207 & -0.030 &  0.184 &  0.025 \\
NVDA  &  0.123 &  0.119 &  0.149 &  0.160 &  0.177 &  0.166 \\
TSLA  &  0.006 &  0.164 &  0.311 &  0.242 &  0.208 &  0.285 \\
NFLX  &  0.098 &  0.034 & -0.103 & -0.162 &  0.152 &  0.043 \\
AMD   &  0.244 &  0.178 &  0.266 &  0.132 &  0.190 &  0.078 \\
GOOGL &  0.015 &  0.198 & -0.253 &  0.043 & -0.071 &  0.156 \\
META  &  0.087 &  0.064 &  0.158 &  0.080 &  0.096 &  0.060 \\
\bottomrule
\end{tabular}
\end{table}

\subsection{Fuzzy Mean--Variance--Skewness Portfolio Problem}
\label{subsec:fuzzy-mvs}

Consider a market of $n$ assets. The rate of return for each asset $i$ is modeled by a fuzzy number $\xi_i$, and $w_i\in \mathbb{R}$ are decision variables representing the percentage share of the decision maker's overall capital allocated to asset $i.$
Assuming $ \xi_1, \xi_2, \ldots, \xi_n$ to be fuzzy numbers, the aggregate return of the portfolio is also given by a fuzzy number $\bigoplus_{i=1}^n w_{i}\xi_{i}$.

Like in the probabilistic model, the mean value $\E\big[\bigoplus_{i=1}^n w_{i}\xi_{i}\big]$ models the expected return for the entire portfolio, while the aggregated risk is characterized by $\V\big[\bigoplus_{i=1}^n w_{i}\xi_{i}\big]$, and aggregated skewness by $\Skew\big[\bigoplus_{i=1}^n w_{i}\xi_{i}\big]$.

A caveat applies to this aggregation: the level-set arithmetic $\oplus$ adds marginal quantiles level by level, implicitly treating the assets as comonotonic and ignoring their empirical dependence structure. This convention is standard in the fuzzy portfolio literature, but should be kept in mind below, where the selected portfolios are evaluated probabilistically.

We consider the following optimization problem:
\begin{equation}\label{eq:PortfolioOptimization}
	\begin{aligned}
		  & \max \; \Skew\big[w_1 \xi_1 \oplus \ldots \oplus w_n \xi_n\big] \\
		  & \text{subject to}                                                             \\
		  & \E\big[w_1 \xi_1 \oplus \ldots \oplus w_n \xi_n\big] \geq \rho_{0}      \\
		  & \V\big[w_1 \xi_1 \oplus \ldots \oplus w_n \xi_n\big] \leq \beta_{0}       \\
		  & w_1 + w_2 + \ldots + w_n = 1                                   \\
		  & 0 \leq w_i \leq 1, \; i = 1, \ldots, n,                                       
	\end{aligned}
\end{equation}
where $\oplus$ denotes fuzzy addition in the level-set sense described in Section~\ref{sec:preliminaries}. 

To establish baseline thresholds, we compute the fuzzy mean and variance of the equally weighted portfolio, $w^{eq} = (1/n,\ldots,1/n)^\top$. Denoting its fuzzy return as $\Xi_{w^{eq}}$, the baseline metrics are defined via the possibilistic mean \eqref{eq:BermMean} and the squared downside risk \eqref{eq:VarBermudez},
\begin{equation}\label{eq:BaselineMeasures}
\begin{aligned}
\bar{\mu} &= \E(\Xi_{w^{eq}}) = \frac{1}{2} \int_{0}^{1} \Big[ (\Xi_{w^{eq}})_d(\alpha) + (\Xi_{w^{eq}})_u(\alpha) \Big] d\alpha, \\
\bar{\sigma}^2 &= \V(\Xi_{w^{eq}}) = \left( \int_{0}^{1} \Big[ (\Xi_{w^{eq}})_u(\alpha) - (\Xi_{w^{eq}})_d(\alpha) \Big] d\alpha \right)^2.
\end{aligned}
\end{equation}
We then set $\rho_{0} = \rho\, \bar{\mu}$ and $\beta_{0} = \beta\, \bar{\sigma}^2$, so that $(\rho,\beta)=(1,1)$ corresponds to imposing at least the baseline fuzzy mean and at most the baseline fuzzy variance. Values $\rho<1$ and $\beta>1$ relax the mean and risk constraints, respectively, while $\rho>1$ and $\beta<1$ tighten them.


\subsection{Optimization Setup and Probabilistic Evaluation}
\label{subsec:setup}

For each $(\rho,\beta)$ pair in a grid $\rho \in \{0.95,1.00,1.05\}$ and $\beta \in \{0.95,1.00,1.05\}$, we approximate the feasible set of portfolios by a large Monte--Carlo ``universe'' of randomly generated weights ($S=50{,}000$ drawn  uniformly on the $n$-asset simplex). For each such portfolio, we compute $\E(\Xi_w)$, $\V(\Xi_w)$, and $\Skew(\Xi_w)$ under each skewness definition. The mean and variance operators in the constraints of \eqref{eq:PortfolioOptimization} are framework-specific: $\Skew_{VB13}$ and all $\Skew_{JKPT}$ configurations use the possibilistic mean \eqref{eq:BermMean} and squared downside risk \eqref{eq:VarBermudez} (the $\Skew_{JKPT}$ family defines skewness only and adopts these operators), whereas $\Skew_{LGY15}$ portfolios use the level-weighted operators \eqref{eq:LiMean} and \eqref{eq:LiVar} of \cite{li2015skewness}. The feasible sets therefore differ: in Table~\ref{tab:results_095_100}, all $\Skew_{VB13}$ and $\Skew_{JKPT}$ rows share the feasible count $12{,}342$ of the $S=50{,}000$ draws, while the $\Skew_{LGY15}$ row differs ($15{,}175$).

Naturally, this random-search scheme approximates, rather than solves, problem \eqref{eq:PortfolioOptimization}: $50{,}000$ draws cover the nine-asset simplex only coarsely, and identical allocations across adjacent configurations in Table~\ref{tab:results_095_100} are a symptom of this discreteness. A local refinement step or a dedicated solver is a natural extension.

Among the feasible portfolios satisfying the fuzzy constraints, we select the one with maximal fuzzy skewness. For $\Skew_{JKPT_1}$ and $\Skew_{JKPT_2}$ we consider a grid of convex combinations $v \in \{0.0, 0.1, \ldots, 1.0\}$ trading off the inner and outer skewness components, with $\alpha \in \{0.10,0.25\}$ for $\Skew_{JKPT_1}$. 

To assess the quality of the selected fuzzy portfolios in standard probabilistic terms, we split the sample into a training period (first two thirds) and an out-of-sample period (remaining third). The chosen portfolios $w^*_{\Skew,\rho,\beta}$ are evaluated out-of-sample on test returns $R_w^{\text{test}}$ by computing mean return $\mu$, standard deviation $\sigma$, Fisher--Pearson probabilistic skewness, and the Sharpe ratio $\mu/\sigma$ (reporting units are stated in the caption of Table~\ref{tab:results_095_100}).

\subsection{Results and Discussion}
\label{subsec:results}

Table~\ref{tab:results_095_100} reports out-of-sample allocations and performance for the representative constraint configuration $(\rho,\beta)=(0.95, 1.00)$; the corresponding Sharpe-ratio profile across all individual configurations is plotted in the Supplementary Material. Several observations follow.

First, portfolios selected under the proposed quantile-based coefficients are competitive with the moment-based benchmarks: the best configuration, $\Skew_{JKPT_1}(0.10)$ with $v \in [0.4, 0.7]$, attains an out-of-sample Sharpe ratio of $0.361$, and the $\Skew_{JKPT_1}(0.25)$ plateau for $v \leq 0.3$ attains $0.360$, versus $0.279$ for $\Skew_{VB13}$ and $0.324$ for $\Skew_{LGY15}$. These figures, however, span a grid of $(\alpha, v)$ configurations inspected ex post: they show that the parameter family contains well-performing members, not that a single, prespecified configuration dominates the benchmarks. A preregistered choice of $(\alpha, v)$, or selection on a validation subperiod followed by a one-shot out-of-sample evaluation, is left to future work.

Second, the out-of-sample window comprises 32 monthly observations from a single training/test split. Sharpe ratio differences of the observed magnitude lie within sampling variability at this sample size, and we have not performed formal significance tests. The results are therefore an illustration of the methodology rather than confirmatory evidence of outperformance.

Third, Table~\ref{tab:asset_statistics_competitors} shows that the quantile-based construction of Section~\ref{subsec:data} drives $\Skew_{LGY15}$ to values of order $10^{-4}$ and below, i.e., numerical noise. Maximizing such an objective selects an essentially arbitrary feasible portfolio, so the comparison documents that $\Skew_{LGY15}$ degenerates for fuzzy numbers constructed from empirical quantiles, rather than evidence of outperforming a meaningful competitor.

Fourth, the selected portfolios are strongly concentrated, up to $62\%$ in a single asset, within a universe of nine growth equities chosen, with hindsight, for their pronounced asymmetry during a predominantly bullish period; the reported performance is specific to this environment. Robustness checks on defensive or cross-asset universes, alternative sample periods, and rolling windows are deferred to a follow-up study.

Finally, high Sharpe ratios concentrate at low-to-moderate values of $v$ ($v \leq 0.7$ for $\Skew_{JKPT_1}(0.10)$, $v \leq 0.3$ for $\Skew_{JKPT_1}(0.25)$, and $v = 0$ for $\Skew_{JKPT_2}$), suggesting that the inner (component-level) skewness carries useful information; the evidence does not support the stronger claim that smaller values of $v$ are uniformly preferable.

\begin{table}[htbp]
\centering
\caption{Out-of-sample portfolio performance and asset allocations for the constraint configuration $(\rho,\beta)=(0.95, 1.00)$: percentage allocations ($w_i$) for the nine assets, and the portfolio's out-of-sample probabilistic mean return ($\mu$, percent per month), Fisher--Pearson skewness, and monthly non-annualized Sharpe ratio (zero risk-free rate). ``Feas.'' is the number of the $S=50{,}000$ random portfolios satisfying the framework-specific constraints (Section~\ref{subsec:setup}). Color highlights mark $\Skew_{JKPT}$ configurations strictly outperforming the better of the two baselines ($\Skew_{VB13}$, $\Skew_{LGY15}$) in a given metric (no configuration strictly exceeds the competitors' best mean, hence no mean cells are highlighted). Rows with identical allocations across adjacent values of $v$ are grouped; the full per-configuration table is provided in the Supplementary Material, and full-grid results with replication code in the companion repository.}
\label{tab:results_095_100}
\footnotesize
\setlength{\tabcolsep}{4pt}

\textbf{Legend:} \colorbox{blue!20}{Beats Comp. Skew} | \colorbox{yellow!40}{Beats Comp. Sharpe}
\vspace{0.5em}

\begin{adjustbox}{max width=\textwidth}
\begin{tabular}{l r rrrrrrr rrrrr}
\toprule
Method & Feas. & AAPL & MSFT & AMZN & NVDA & TSLA & NFLX & AMD & GOOGL & META & Mean & Skew & Sharpe \\ \midrule
VB13 & 12342 & 0.01 & 0.03 & 0.00 & 0.02 & 0.27 & 0.07 & 0.01 & 0.35 & 0.23 & 1.695 & -0.225 & 0.279 \\
LGY15 & 15175 & 0.01 & 0.04 & 0.01 & 0.62 & 0.08 & 0.03 & 0.15 & 0.06 & 0.00 & 3.498 & -0.287 & 0.324 \\
JKPT1\_{0.10} $v \in [0.0, 0.2]$ & 12342 & 0.04 & 0.06 & 0.28 & 0.33 & 0.00 & 0.00 & 0.22 & 0.00 & 0.08 & 3.415 & -0.515 & \cellcolor{yellow!40}0.358 \\
JKPT1\_{0.10} $v = 0.3$ & 12342 & 0.00 & 0.03 & 0.16 & 0.39 & 0.01 & 0.03 & 0.19 & 0.08 & 0.11 & 3.264 & -0.461 & \cellcolor{yellow!40}0.340 \\
JKPT1\_{0.10} $v \in [0.4, 0.7]$ & 12342 & 0.00 & 0.02 & 0.23 & 0.28 & 0.07 & 0.02 & 0.20 & 0.17 & 0.01 & 3.137 & -0.364 & \cellcolor{yellow!40}0.361 \\
JKPT1\_{0.10} $v \in [0.8, 1.0]$ & 12342 & 0.01 & 0.00 & 0.02 & 0.01 & 0.27 & 0.01 & 0.08 & 0.56 & 0.05 & 1.813 & \cellcolor{blue!20}-0.198 & 0.303 \\
JKPT1\_{0.25} $v \in [0.0, 0.3]$ & 12342 & 0.01 & 0.04 & 0.48 & 0.11 & 0.15 & 0.00 & 0.09 & 0.00 & 0.10 & 2.664 & \cellcolor{blue!20}-0.109 & \cellcolor{yellow!40}0.360 \\
JKPT1\_{0.25} $v \in [0.4, 0.6]$ & 12342 & 0.00 & 0.06 & 0.17 & 0.38 & 0.26 & 0.00 & 0.03 & 0.03 & 0.07 & 2.722 & \cellcolor{blue!20}-0.068 & 0.320 \\
JKPT1\_{0.25} $v \in [0.7, 0.8]$ & 12342 & 0.04 & 0.05 & 0.05 & 0.53 & 0.21 & 0.01 & 0.05 & 0.03 & 0.03 & 3.032 & \cellcolor{blue!20}-0.103 & 0.312 \\
JKPT1\_{0.25} $v = 0.9$ & 12342 & 0.02 & 0.05 & 0.02 & 0.57 & 0.21 & 0.01 & 0.00 & 0.13 & 0.00 & 2.821 & \cellcolor{blue!20}-0.022 & 0.289 \\
JKPT1\_{0.25} $v = 1.0$ & 12342 & 0.04 & 0.02 & 0.02 & 0.13 & 0.35 & 0.00 & 0.00 & 0.37 & 0.06 & 1.892 & -0.228 & 0.273 \\
JKPT2 $v = 0.0$ & 12342 & 0.00 & 0.02 & 0.62 & 0.13 & 0.07 & 0.06 & 0.06 & 0.02 & 0.01 & 2.777 & \cellcolor{blue!20}-0.063 & \cellcolor{yellow!40}0.347 \\
JKPT2 $v \in [0.1, 0.4]$ & 12342 & 0.00 & 0.01 & 0.23 & 0.51 & 0.17 & 0.00 & 0.02 & 0.02 & 0.03 & 2.937 & \cellcolor{blue!20}-0.059 & 0.305 \\
JKPT2 $v \in [0.5, 0.7]$ & 12342 & 0.01 & 0.06 & 0.02 & 0.35 & 0.32 & 0.04 & 0.00 & 0.18 & 0.03 & 2.452 & \cellcolor{blue!20}-0.032 & 0.295 \\
JKPT2 $v \in [0.8, 1.0]$ & 12342 & 0.04 & 0.02 & 0.02 & 0.13 & 0.35 & 0.00 & 0.00 & 0.37 & 0.06 & 1.892 & -0.228 & 0.273 \\
\bottomrule
\end{tabular}
\end{adjustbox}
\end{table}


\subsection{Analytical Time Complexity Analysis}
\label{subsec:complexity}
To evaluate the skewness of a fuzzy portfolio, the optimization algorithm must first construct the aggregate portfolio fuzzy number, defined as $\tilde{\xi}_P = \bigoplus_{i} w_i \tilde{\xi}_i$. Subsequently, the corresponding skewness coefficient is computed. To evaluate the exact computational efficiency of this latter step, let $N$ denote the number of discrete $\alpha$-levels utilized for numerical integration (e.g., $N=21$ for an integration step $\Delta\alpha=0.05$, or $N=101$ for $\Delta\alpha=0.01$).

The competing indicators $\Skew_{VB13}$ and $\Skew_{LGY15}$ rely on a moment-based framework. Computing the third central fuzzy moment requires a two-stage numerical procedure:
\begin{enumerate}
    \item Evaluation of the fuzzy expected value $\E$, which necessitates numerical integration across the $N$ discrete $\alpha$-levels.
    \item Integration of the cubed deviations, $\int (x - \E)^3$. For each of the $N$ discrete levels, the algorithm must subtract the previously computed expected value $\E$, cube the resulting difference, scale it by the appropriate integration weight, and accumulate the sum.
\end{enumerate}
Consequently, the computational cost scales linearly with the resolution of the fuzzy number, yielding an analytical time complexity of $\mathcal{O}(N)$.

The proposed integral-based coefficient $\Skew_{JKPT2}$ similarly necessitates a numerical integration procedure over the continuous domain of $\alpha$-levels. As a result, its computational expenditure depends directly on the grid resolution, leading to a general time complexity of $\mathcal{O}(N)$.

In contrast, the proposed point-based coefficient $\Skew_{JKPT1}$ circumvents the numerical integration loop entirely. Its formulation requires evaluating the membership function at a strictly bounded, finite set of points, typically three or four distinct levels (e.g., $\alpha \in \{0.25, 0.50, 0.75, 1.00\}$).

The computational execution reduces to a direct lookup of these few scalar values followed by a deterministic, fixed sequence of basic arithmetic operations (additions and divisions). Since the execution time is completely decoupled from the discretization resolution $N$, the analytical time complexity of $\Skew_{JKPT1}$ is strictly $\mathcal{O}(1)$ (constant time).

We emphasize the scope of this analysis: it concerns a single skewness evaluation. In the experiment above, total runtime is dominated by the size of the Monte--Carlo universe of candidate portfolios (and the construction of each aggregate fuzzy number), not by the per-evaluation cost. The $\mathcal{O}(1)$ property of $\Skew_{JKPT_1}$ therefore does not, by itself, establish end-to-end scalability to large universes or real-time use; a measured scaling study, with dedicated optimizers in place of random search, would be required to substantiate such claims and is left for future work.

\color{black}

\section{Conclusion}\label{sec:conclusions}
This article introduced a nonmoment-based and quantile-based skewness measure for fuzzy numbers, addressing an open problem in the existing literature and complementing the established moment-based approaches. We established a probabilistic foundation for fuzzy numbers by interpreting their left and right membership function components as cumulative and survival probability functions, respectively. This approach allows for a direct interpretation of $\alpha$-cuts as generalized quantiles, moving beyond the prevalent, but often theoretically unfounded, practice of directly substituting fuzzy membership functions for probability density functions in moment-based skewness formulas.

The proposed $\Skew_{{JKPT}_1}$ and $\Skew_{{JKPT}_2}$ coefficients combine two complementary constituents: an ``inner'' component quantifying the intrinsic asymmetry of the underlying probabilistic generative processes ($X_L$ and $X_R$), and an ``outer'' component capturing the emergent skewness of the fuzzy number's overall profile. This dual perspective, combined with a flexible weighting factor $v \in [0,1]$, offers additional insight into the asymmetry of fuzzy quantities.

\noindent We demonstrated that the new coefficients possess several desirable mathematical properties, crucial for robust application. These include a bounded range of $[-1, 1]$, invariance under both scale and location transformations, symmetry (yielding zero for perfectly symmetric fuzzy numbers), and monotonicity (correctly identifying the direction of skewness). Furthermore, their quantile-based nature renders them inherently robust with respect to outliers, a common challenge for moment-based coefficients.


The practical use of our proposed skewness measures was illustrated within a fuzzy mean--variance--skewness portfolio optimization framework using real financial data. By evaluating a sizable universe of feasible portfolios generated via a Monte--Carlo approach under fuzzy constraints, our method consistently provided sensible asset allocations. The analytical time complexity analysis shows that while the moment-based competitors ($\Skew_{LGY15}$~\cite{li2015skewness} and $\Skew_{VB13}$~\cite{VB13}) and our integral-based $\Skew_{JKPT_2}$ variant scale linearly with the discretization resolution ($\mathcal{O}(N)$), the point-based $\Skew_{JKPT_1}$ variant operates in constant time ($\mathcal{O}(1)$). This decoupling from numerical integration loops reduces the per-evaluation cost. Furthermore, our empirical investigation exposed a practical limitation of the $\Skew_{LGY15}$ coefficient: due to its rigid cubic deviation from the mean, its values collapse near zero, limiting its effectiveness for fuzzy numbers constructed via empirical quantiles. In contrast, the proposed $\Skew_{JKPT}$ indices captured the tail asymmetry of the data and achieved competitive out-of-sample risk-adjusted performance in our illustration, subject to the caveats discussed in Section~\ref{subsec:results}.

This work provides an interpretable and computationally efficient tool for decision-makers operating under uncertainty, particularly in financial contexts where higher order risk preferences matter. It illustrates how a probabilistically grounded approach to fuzzy skewness can inform practical decision-making in financial engineering.

Future research could explore several directions. For instance, optimal selection of weighting factor $v$ and quantile levels $\alpha$ for $\Skew_{JKPT_1}(\alpha)(v)$ across application scenarios (e.g., sensitivity analyses and machine learning). Applying these coefficients to other fuzzy decision-making domains (e.g. supply chain risk and medical diagnosis). Analytical comparison of $\Skew_{JKPT_1}$ and $\Skew_{JKPT_2}$ behavior for specific fuzzy number classes (e.g., from different probability distributions or histograms). Developing a complementary quantile-based kurtosis measure to provide a more complete description of the underlying fuzzy distributions. Examining the coefficients' behavior for parameterized families of fuzzy numbers, such as triangular dense fuzzy sets \cite{De2016469} and dense fuzzy lock sets \cite{De20187243}, where designated key parameters (``locks'') within the set allow the decision maker to dynamically control or reset the system. Extending the proposed framework to multivariable fuzzy sets, e.g., via joint quantile functions, to capture dependencies in higher dimensions. Applying the framework to larger and more diverse asset universes, with rolling-window evaluations and formal significance tests, and integrating the coefficients into dedicated optimization algorithms.

Finally, further theoretical exploration could focus on the strict separation between skewness and the definitions of mean and variance, ensuring a more robust mathematical foundation for shape analysis in fuzzy environments.

\bibliography{Quantile_skewness}
\bibliographystyle{unsrturl}

\clearpage
\appendix
\renewcommand{\thesection}{S.\Roman{section}}

\section*{Supplementary Material}
\addcontentsline{toc}{section}{Supplementary Material}

This appendix contains the supplementary material accompanying the article: (\ref{sec:supp-example}) a step-by-step visualization example of the construction of a fuzzy number from two probability distributions, (\ref{sec:supp-table}) the full, per-configuration out-of-sample results table for the constraint configuration $(\rho,\beta)=(0.95,1.00)$, of which the main text presents a condensed version, together with the corresponding Sharpe-ratio profile plot, and (\ref{sec:supp-asymptotics}) an asymptotic analysis of the ratios between the moment-based benchmark indices and the proposed quantile-based indices under increasing scale and one-sided spread.

\section{Visualizing the construction of a fuzzy number from two probability distributions}
\label{sec:supp-example}

The following example (referred to in Section~\ref{subsec:JKPT_Skewness_Examples} of the main text) visualizes the constructive procedure of the paper, transforming two independent random variables into a single fuzzy number using relatively small, manageable numbers.

\begin{example}{\textit{Visualizing the construction of a fuzzy number from two probability distributions}.}\label{ex:RV_to_FN_Visualization}
In this scenario, we construct a fuzzy number $\xi$ representing a moderate positive skew. We utilize two Beta distributions, $X_L$ and $X_R$, defined on adjacent supports $[100, 102]$ and $[102, 110]$ respectively. Both share the same shape parameters ($\alpha=0.5, \beta=2$), creating a consistent generative logic, but the difference in support width induces the asymmetry.
The calculated coefficients for this fuzzy number are: $\Skew_{JKPT_{1}}(0.1) (0.5) = 0.461,$ $\Skew_{JKPT_{1}}(0.25) (0.5) = 0.235,$  $\Skew_{JKPT_{2}}(0.5) = 0.347,$ $\Skew_{VB13} = 0.298,$  $\Skew_{LGY15} = 1.96$.

\smallskip
In this ``small number'' regime, all coefficients accordantly identify positive skewness. $\Skew_{JKPT}$ and $\Skew_{VB13}$ provide values that are comparable in magnitude, falling within a standard interpretable range. $\Skew_{LGY15}$ is already notably higher, but not yet extreme.
\end{example}

\begin{figure}[ht]
  \centering
  \includegraphics[width=0.7\textwidth]{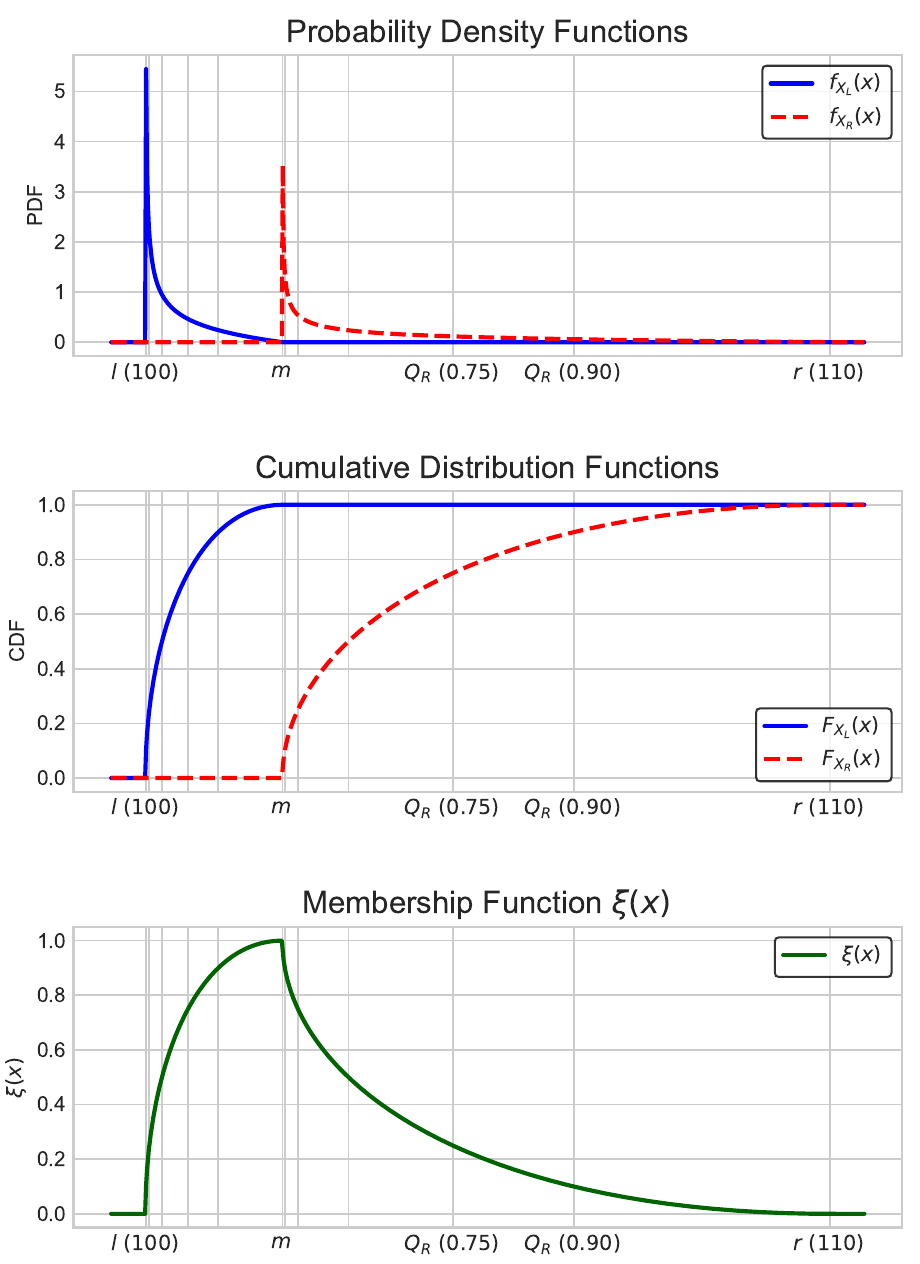}
  \caption{Visualization for Example~\ref{ex:RV_to_FN_Visualization}.
\newline\textbf{Top:} PDFs of the generative Beta distributions.
\newline\textbf{Middle:} Corresponding CDFs.
\newline\textbf{Bottom:} The resulting membership function $\xi(x)$.
\newline The x-axis tick marks indicate the core/support boundaries and the quantile points used for the point-based skewness calculation.}
  \label{Fig:Ex1_Visual}
\end{figure}

\section{Full out-of-sample results table and Sharpe-ratio profile}
\label{sec:supp-table}

\begin{figure}[htbp]
  \centering
  \includegraphics[width=0.8\textwidth]{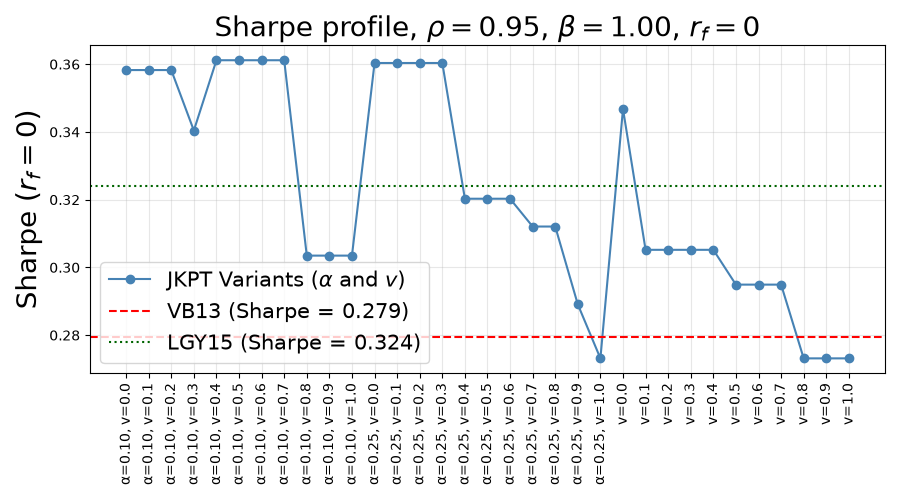}
  \caption{Out-of-sample Sharpe Ratio profile for portfolios optimized under various $\Skew_{JKPT}$ configurations for the constraint setting $(\rho,\beta)=(0.95, 1.00)$. The x-axis enumerates specific $\Skew_{JKPT_1}$ (for $\alpha=0.10, 0.25$) and $\Skew_{JKPT_2}$ variants across their respective $v$ parameter ranges. The y-axis represents the achieved out-of-sample Sharpe Ratio. Horizontal lines indicate the Sharpe Ratios obtained by portfolios optimized using the $\Skew_{VB13}$ and $\Skew_{LGY15}$ skewness definitions under the same fuzzy constraints. The profile shows how out-of-sample performance varies with the configuration $(\alpha, v)$: several configurations attain higher Sharpe ratios than the benchmarks, while others do not; see Section~\ref{sec:empirical} of the main text for a discussion.}
  \label{Fig:sharpe-profile}
\end{figure}

Table~\ref{tab:results_full} lists the out-of-sample portfolio performance and asset allocations for every individual configuration $(\alpha, v)$ under the constraint configuration $(\rho,\beta)=(0.95,1.00)$; the main text (Section~\ref{sec:empirical}) presents a condensed version where rows with identical allocations across consecutive values of $v$ are aggregated. Mean returns are in percent per month; Sharpe ratios are monthly and non-annualized, with a zero risk-free rate. ``Feas.'' is the number of the $S=50{,}000$ random portfolios satisfying the framework-specific constraints. Color highlights denote instances where a $\Skew_{JKPT}$ configuration strictly outperforms the maximum value achieved by the baseline competitors ($\Skew_{VB13}$ and $\Skew_{LGY15}$) in a given metric (Skew or Sharpe). Fig.~\ref{Fig:sharpe-profile} plots the corresponding Sharpe-ratio profile.

\begin{table}[htbp]
\centering
\caption{Out-of-sample portfolio performance and asset allocations for the constraint configuration $(\rho,\beta)=(0.95, 1.00)$: full per-configuration results.}
\label{tab:results_full}
\footnotesize
\setlength{\tabcolsep}{4pt}

\textbf{Legend:} \colorbox{blue!20}{Beats Comp. Skew} | \colorbox{yellow!40}{Beats Comp. Sharpe}
\vspace{0.5em}

\begin{adjustbox}{max width=\textwidth}
\begin{tabular}{l r rrrrrrr rrrrr}
\toprule
Method & Feas. & AAPL & MSFT & AMZN & NVDA & TSLA & NFLX & AMD & GOOGL & META & Mean & Skew & Sharpe \\ \midrule
VB13 & 12342 & 0.01 & 0.03 & 0.00 & 0.02 & 0.27 & 0.07 & 0.01 & 0.35 & 0.23 & 1.695 & -0.225 & 0.279 \\
LGY15 & 15175 & 0.01 & 0.04 & 0.01 & 0.62 & 0.08 & 0.03 & 0.15 & 0.06 & 0.00 & 3.498 & -0.287 & 0.324 \\
JKPT1\_{0.10} v=0.0 & 12342 & 0.04 & 0.06 & 0.28 & 0.33 & 0.00 & 0.00 & 0.22 & 0.00 & 0.08 & 3.415 & -0.515 & \cellcolor{yellow!40}0.358 \\
JKPT1\_{0.10} v=0.1 & 12342 & 0.04 & 0.06 & 0.28 & 0.33 & 0.00 & 0.00 & 0.22 & 0.00 & 0.08 & 3.415 & -0.515 & \cellcolor{yellow!40}0.358 \\
JKPT1\_{0.10} v=0.2 & 12342 & 0.04 & 0.06 & 0.28 & 0.33 & 0.00 & 0.00 & 0.22 & 0.00 & 0.08 & 3.415 & -0.515 & \cellcolor{yellow!40}0.358 \\
JKPT1\_{0.10} v=0.3 & 12342 & 0.00 & 0.03 & 0.16 & 0.39 & 0.01 & 0.03 & 0.19 & 0.08 & 0.11 & 3.264 & -0.461 & \cellcolor{yellow!40}0.340 \\
JKPT1\_{0.10} v=0.4 & 12342 & 0.00 & 0.02 & 0.23 & 0.28 & 0.07 & 0.02 & 0.20 & 0.17 & 0.01 & 3.137 & -0.364 & \cellcolor{yellow!40}0.361 \\
JKPT1\_{0.10} v=0.5 & 12342 & 0.00 & 0.02 & 0.23 & 0.28 & 0.07 & 0.02 & 0.20 & 0.17 & 0.01 & 3.137 & -0.364 & \cellcolor{yellow!40}0.361 \\
JKPT1\_{0.10} v=0.6 & 12342 & 0.00 & 0.02 & 0.23 & 0.28 & 0.07 & 0.02 & 0.20 & 0.17 & 0.01 & 3.137 & -0.364 & \cellcolor{yellow!40}0.361 \\
JKPT1\_{0.10} v=0.7 & 12342 & 0.00 & 0.02 & 0.23 & 0.28 & 0.07 & 0.02 & 0.20 & 0.17 & 0.01 & 3.137 & -0.364 & \cellcolor{yellow!40}0.361 \\
JKPT1\_{0.10} v=0.8 & 12342 & 0.01 & 0.00 & 0.02 & 0.01 & 0.27 & 0.01 & 0.08 & 0.56 & 0.05 & 1.813 & \cellcolor{blue!20}-0.198 & 0.303 \\
JKPT1\_{0.10} v=0.9 & 12342 & 0.01 & 0.00 & 0.02 & 0.01 & 0.27 & 0.01 & 0.08 & 0.56 & 0.05 & 1.813 & \cellcolor{blue!20}-0.198 & 0.303 \\
JKPT1\_{0.10} v=1.0 & 12342 & 0.01 & 0.00 & 0.02 & 0.01 & 0.27 & 0.01 & 0.08 & 0.56 & 0.05 & 1.813 & \cellcolor{blue!20}-0.198 & 0.303 \\
JKPT1\_{0.25} v=0.0 & 12342 & 0.01 & 0.04 & 0.48 & 0.11 & 0.15 & 0.00 & 0.09 & 0.00 & 0.10 & 2.664 & \cellcolor{blue!20}-0.109 & \cellcolor{yellow!40}0.360 \\
JKPT1\_{0.25} v=0.1 & 12342 & 0.01 & 0.04 & 0.48 & 0.11 & 0.15 & 0.00 & 0.09 & 0.00 & 0.10 & 2.664 & \cellcolor{blue!20}-0.109 & \cellcolor{yellow!40}0.360 \\
JKPT1\_{0.25} v=0.2 & 12342 & 0.01 & 0.04 & 0.48 & 0.11 & 0.15 & 0.00 & 0.09 & 0.00 & 0.10 & 2.664 & \cellcolor{blue!20}-0.109 & \cellcolor{yellow!40}0.360 \\
JKPT1\_{0.25} v=0.3 & 12342 & 0.01 & 0.04 & 0.48 & 0.11 & 0.15 & 0.00 & 0.09 & 0.00 & 0.10 & 2.664 & \cellcolor{blue!20}-0.109 & \cellcolor{yellow!40}0.360 \\
JKPT1\_{0.25} v=0.4 & 12342 & 0.00 & 0.06 & 0.17 & 0.38 & 0.26 & 0.00 & 0.03 & 0.03 & 0.07 & 2.722 & \cellcolor{blue!20}-0.068 & 0.320 \\
JKPT1\_{0.25} v=0.5 & 12342 & 0.00 & 0.06 & 0.17 & 0.38 & 0.26 & 0.00 & 0.03 & 0.03 & 0.07 & 2.722 & \cellcolor{blue!20}-0.068 & 0.320 \\
JKPT1\_{0.25} v=0.6 & 12342 & 0.00 & 0.06 & 0.17 & 0.38 & 0.26 & 0.00 & 0.03 & 0.03 & 0.07 & 2.722 & \cellcolor{blue!20}-0.068 & 0.320 \\
JKPT1\_{0.25} v=0.7 & 12342 & 0.04 & 0.05 & 0.05 & 0.53 & 0.21 & 0.01 & 0.05 & 0.03 & 0.03 & 3.032 & \cellcolor{blue!20}-0.103 & 0.312 \\
JKPT1\_{0.25} v=0.8 & 12342 & 0.04 & 0.05 & 0.05 & 0.53 & 0.21 & 0.01 & 0.05 & 0.03 & 0.03 & 3.032 & \cellcolor{blue!20}-0.103 & 0.312 \\
JKPT1\_{0.25} v=0.9 & 12342 & 0.02 & 0.05 & 0.02 & 0.57 & 0.21 & 0.01 & 0.00 & 0.13 & 0.00 & 2.821 & \cellcolor{blue!20}-0.022 & 0.289 \\
JKPT1\_{0.25} v=1.0 & 12342 & 0.04 & 0.02 & 0.02 & 0.13 & 0.35 & 0.00 & 0.00 & 0.37 & 0.06 & 1.892 & -0.228 & 0.273 \\
JKPT2 v=0.0 & 12342 & 0.00 & 0.02 & 0.62 & 0.13 & 0.07 & 0.06 & 0.06 & 0.02 & 0.01 & 2.777 & \cellcolor{blue!20}-0.063 & \cellcolor{yellow!40}0.347 \\
JKPT2 v=0.1 & 12342 & 0.00 & 0.01 & 0.23 & 0.51 & 0.17 & 0.00 & 0.02 & 0.02 & 0.03 & 2.937 & \cellcolor{blue!20}-0.059 & 0.305 \\
JKPT2 v=0.2 & 12342 & 0.00 & 0.01 & 0.23 & 0.51 & 0.17 & 0.00 & 0.02 & 0.02 & 0.03 & 2.937 & \cellcolor{blue!20}-0.059 & 0.305 \\
JKPT2 v=0.3 & 12342 & 0.00 & 0.01 & 0.23 & 0.51 & 0.17 & 0.00 & 0.02 & 0.02 & 0.03 & 2.937 & \cellcolor{blue!20}-0.059 & 0.305 \\
JKPT2 v=0.4 & 12342 & 0.00 & 0.01 & 0.23 & 0.51 & 0.17 & 0.00 & 0.02 & 0.02 & 0.03 & 2.937 & \cellcolor{blue!20}-0.059 & 0.305 \\
JKPT2 v=0.5 & 12342 & 0.01 & 0.06 & 0.02 & 0.35 & 0.32 & 0.04 & 0.00 & 0.18 & 0.03 & 2.452 & \cellcolor{blue!20}-0.032 & 0.295 \\
JKPT2 v=0.6 & 12342 & 0.01 & 0.06 & 0.02 & 0.35 & 0.32 & 0.04 & 0.00 & 0.18 & 0.03 & 2.452 & \cellcolor{blue!20}-0.032 & 0.295 \\
JKPT2 v=0.7 & 12342 & 0.01 & 0.06 & 0.02 & 0.35 & 0.32 & 0.04 & 0.00 & 0.18 & 0.03 & 2.452 & \cellcolor{blue!20}-0.032 & 0.295 \\
JKPT2 v=0.8 & 12342 & 0.04 & 0.02 & 0.02 & 0.13 & 0.35 & 0.00 & 0.00 & 0.37 & 0.06 & 1.892 & -0.228 & 0.273 \\
JKPT2 v=0.9 & 12342 & 0.04 & 0.02 & 0.02 & 0.13 & 0.35 & 0.00 & 0.00 & 0.37 & 0.06 & 1.892 & -0.228 & 0.273 \\
JKPT2 v=1.0 & 12342 & 0.04 & 0.02 & 0.02 & 0.13 & 0.35 & 0.00 & 0.00 & 0.37 & 0.06 & 1.892 & -0.228 & 0.273 \\
\bottomrule
\end{tabular}
\end{adjustbox}
\end{table}

\section{Asymptotic behavior of the skewness indices}
\label{sec:supp-asymptotics}

This section formalizes the behavior illustrated by the stress-test example of the main text (Example~\ref{ex:RV_to_FN_ExtremeSkew}, Section~\ref{subsec:JKPT_Skewness_Examples}): we compare the proposed coefficients $\Skew_{JKPT_1}$ and $\Skew_{JKPT_2}$ with the moment-based benchmarks $\Skew_{VB13}$ and $\Skew_{LGY15}$ under two asymptotic regimes and report the behavior of the ratios (benchmark index over proposed index).

\subsection{Scaling regime}

\begin{proposition}[Scale asymptotics]\label{prop:scaling}
Let $\xi$ be a fuzzy number and $k>0$. Then
\begin{align*}
\Skew_{JKPT_1}(k\xi) &= \Skew_{JKPT_1}(\xi), &
\Skew_{JKPT_2}(k\xi) &= \Skew_{JKPT_2}(\xi), \\
\Skew_{VB13}(k\xi) &= \Skew_{VB13}(\xi), &
\Skew_{LGY15}(k\xi) &= k^3\, \Skew_{LGY15}(\xi).
\end{align*}
Consequently, whenever $\Skew_{LGY15}(\xi) \neq 0$,
\[
\frac{\Skew_{LGY15}(k\xi)}{\Skew_{JKPT_2}(k\xi)}
= k^3 \cdot \frac{\Skew_{LGY15}(\xi)}{\Skew_{JKPT_2}(\xi)}
\;\xrightarrow[k \to \infty]{}\; \infty .
\]
\end{proposition}

\begin{proof}
The level-set representation of $k\xi$ is $[k\xi_d, k\xi_u]$. Numerator and denominator of each inner and outer constituent of $\Skew_{JKPT_1}$ and $\Skew_{JKPT_2}$ are positively homogeneous of degree one in the $\alpha$-cut functions, so the ratios are invariant. For $\Skew_{VB13}$, the third possibilistic moment $\mu_3$ and the cubed downside risk $\omega^3$ are both homogeneous of degree three, so their ratio is invariant. $\Skew_{LGY15}$ cubes deviations from the mean without any normalization and is therefore homogeneous of degree three.
\end{proof}

This is the effect visible in the stress-test example of the main text, where supports of order $10^5$ (instead of order $10^2$ as in Example~\ref{ex:RV_to_FN_Visualization}) drive $\Skew_{LGY15}$ to values of order $10^8$.

\subsection{One-sided spread regime}

We next fix the left side of the fuzzy number and stretch its right support.

\begin{proposition}[One-sided stretch]\label{prop:stretch}
For fixed $\ell > 0$, consider the triangular family $\xi_w$ with $\alpha$-cut functions
$\xi_d(\alpha) = -\ell(1-\alpha)$ and $\xi_u(\alpha) = w(1-\alpha)$, $w \to \infty$. Then:
\begin{enumerate}
\item $\displaystyle \Skew_{JKPT_1}(\xi_w)(\alpha)(v) = \Skew_{JKPT_2}(\xi_w)(v) = v\,\frac{w-\ell}{w+\ell} \;\xrightarrow[w\to\infty]{}\; v \;\le\; 1$;
\item $\displaystyle \Skew_{VB13}(\xi_w) \;\xrightarrow[w\to\infty]{}\; \tfrac{1}{4}$;
\item $\displaystyle \Skew_{LGY15}(\xi_w) = \tfrac{19}{1080}\, w^3 \big(1 + \mathcal{O}(1/w)\big) \;\xrightarrow[w\to\infty]{}\; \infty$.
\end{enumerate}
Hence $\Skew_{LGY15}/\Skew_{JKPT_2} = \Theta(w^3)$, while $\Skew_{VB13}/\Skew_{JKPT_2} \to 1/(4v)$.
\end{proposition}

\begin{proof}
(1) Since the membership functions are linear, the associated random variables are uniformly distributed, $X_L \sim U[-\ell, 0]$ and $X_R \sim U[0, w]$; the Groeneveld--Meeden coefficient of any uniform distribution vanishes (uniform distributions are symmetric), so the inner constituents are zero. For the outer constituents, with core $\{0\}$,
$\int_0^{1/2} [\xi_d(\alpha) + \xi_u(\alpha)]\, d\alpha = \tfrac{3}{8}(w-\ell)$ and
$\int_0^{1/2} [\xi_u(\alpha) - \xi_d(\alpha)]\, d\alpha = \tfrac{3}{8}(w+\ell)$;
the point-based version yields the same ratio at every level $\alpha$.
(2) Here $\E_{VB13} = (w-\ell)/4$ and $\omega = (w+\ell)/2$. Using
$\int_0^1 \big(c(1-\alpha) - E\big)^3 d\alpha = \big((c-E)^4 - E^4\big)/(4c)$ for $c \in \{-\ell, w\}$
and expanding for large $w$ gives $\mu_3 = w^3/32 + \mathcal{O}(w^2)$, while $\omega^3 = w^3/8 + \mathcal{O}(w^2)$.
(3) $\E_{LGY15} = (w-\ell)/6$, and the leading terms of
$\int_0^1 \alpha\big[(\xi_d - E)^3 + (\xi_u - E)^3\big] d\alpha$ are
$w^3 \int_0^1 \alpha\,(5/6-\alpha)^3\, d\alpha = \tfrac{43}{2160} w^3$ and $-\tfrac{5}{2160} w^3$, summing to $\tfrac{19}{1080} w^3$.
\end{proof}

\begin{remark}[Unboundedness of $\Skew_{VB13}$]\label{rem:VB13unbounded}
Although $\Skew_{VB13}$ is scale-invariant (Proposition~\ref{prop:scaling}) and remains bounded along the triangular family of Proposition~\ref{prop:stretch}, it is not uniformly bounded across shapes. Consider fuzzy numbers whose width concentrates near $\alpha = 0$ (a long, thin right tail): take $\xi_d(\alpha) = -\ell_p(1-\alpha)$ and $\xi_u(\alpha) = w(1-\alpha)^p$ with $\ell_p = w/(p+1)$ and $p \to \infty$. Direct computation gives $\E_{VB13} = \tfrac{1}{4}\,w/(p+1)$, $\omega = \tfrac{3}{2}\,w/(p+1)$, and $\mu_3 = \tfrac{1}{6}\,w^3/(p+1)\,(1+o(1))$, so that
\[
\Skew_{VB13} \;=\; \tfrac{4}{81}\,(p+1)^2\,\big(1+o(1)\big) \;\xrightarrow[p\to\infty]{}\; \infty,
\]
while $\Skew_{JKPT_2}(v) \to \tfrac{5}{11}v + \tfrac{1}{2}(1-v) \in (0,1)$ remains bounded (the outer constituent tends to $5/11$; the inner constituent tends to $\tfrac{1}{2}(0+1)$ since $X_L$ is uniform and the Groeneveld--Meeden coefficient of $X_R$ tends to one). This tail-concentration mechanism is precisely the situation of the stress-test example in the main text (right component generated by a $Beta(0.1, 2)$ distribution), where $\Skew_{VB13} = 34.6$.
\end{remark}

In summary, in both regimes the ratio of the benchmark to the proposed index diverges, cubically for $\Skew_{LGY15}$ under scaling and under one-sided stretching, and quadratically in the tail-concentration parameter for $\Skew_{VB13}$, while $\Skew_{JKPT_1}$ and $\Skew_{JKPT_2}$ remain in $[-1, 1]$ with interpretable limits.

\end{document}